\documentclass[12pt]{article}
\usepackage{amsfonts}
\usepackage{graphics}
\newtheorem{THEOREM}{Theorem}[section]
\newtheorem{LEMMA}{Lemma}[section]
\newcommand{\qed}{\ \rule[-1pt]{4pt}{8pt}

                                         \vspace{2ex} }
\newenvironment{PROOF}{

                       \noindent{\bf Proof}.}{\qed}
\newcounter{labelflag} 

\newcommand{\Label}[1]{
                       \ifnum\thelabelflag=1
                          \ifmmode
                             \makebox[0in][l]{\qquad\fbox{\rm#1}}
                          \else
                             \marginpar{\vspace{0.7\baselineskip}
                                        \hspace{-1.1\textwidth}
                                        \fbox{\rm#1}}
                          \fi
                       \fi
                       \label{#1}
                      }
\newcommand{\be}{\begin{equation}}
\newcommand{\ee}{\end{equation}}
\newcommand{\eps}{\varepsilon}
\newcommand{\half}{\textstyle\frac{1}{2}}
\newcommand{\Me}{{\cal M}_\varepsilon}
\newcommand{\Mo}{{\cal M}_0}
\newcommand{\As}{\tilde{A}}
\newcommand{\Bs}{\tilde{B}}
\newcommand{\Us}{\tilde{U}}
\newcommand{\bfLs}{\tilde{\Lambda}}
\newcommand{\psis}{\tilde{\psi}}
\newcommand{\Ks}{\tilde{\mathcal{K}}_\varepsilon}
\newcommand{\Ps}{\tilde{P}}
\newcommand{\Af}{A}
\newcommand{\Bf}{B}
\newcommand{\Uf}{U}
\newcommand{\Lf}{L}
\newcommand{\bfLf}{\Lambda}
\newcommand{\psif}{\psi}
\newcommand{\Kf}{\mathcal{K}_\varepsilon}
\begin{document}

\begin{center}
{\large\textbf{ANALYSIS OF THE CSP REDUCTION METHOD}}\\[1ex]
{\large\textbf{FOR CHEMICAL KINETICS}}\\[1ex]
%Draft -- \today

\vspace{2ex}
Antonios~Zagaris,$^1$
Hans G.~Kaper,$^2$
Tasso J.~Kaper$^1$

$^1$~Department of Mathematics and Center for BioDynamics\\
Boston University,
Boston, Massachusetts, USA \\[1ex]
$^2$~Mathematics and Computer Science Division,\\
Argonne National Laboratory,
Argonne, Illinois, USA \\[1ex]
\end{center}

\begin{small}
\paragraph{Abstract.}
This article is concerned with
the asymptotic accuracy of the
Computational Singular Perturbation (CSP)
method developed by Lam and Goussis
to reduce the dimensionality of a
system of chemical kinetics equations.
The method exploits the presence of
disparate time scales to model
the dynamics by an evolution equation
on a lower-dimensional slow manifold.
In this article it is shown that
the successive applications
of the CSP algorithm generate, 
order by order, 
the asymptotic expansion
of a slow manifold.
The results are illustrated
on the Michaelis--Menten--Henri
equations of enzyme kinetics.



\end{small}
%%%%%%%%%%%%%%%%%%%%%%%%%%%%%%%%%%%%%%%%%%%%%%%%%%%%%%%%
%  Body
%%%%%%%%%%%%%%%%%%%%%%%%%%%%%%%%%%%%%%%%%%%%%%%%%%%%%%%%

\section{Introduction and Summary of Results \label{s-intro}}
\setcounter{equation}{0}
Reduction methods decrease the size
and complexity of systems of kinetic equations.
They are effective
when a small number of variables
can be singled out as evolving
on a ``slow manifold''
and the remaining (fast) variables
somehow follow from the slow variables.
In such cases, the system of kinetic equations
can be reduced to a much smaller system
for the evolution of only the slow variables,
and the fast variables can be determined simply
by table look-ups or by direct computation.
Over the years, a large number of reduction methods
have been proposed and implemented in computer codes;
references can be found
in our earlier article~\cite{KK-2002},
and additional references
are~\cite{CFL-2002, GK-2002, M-2002}.

The focus of Ref.~\cite{KK-2002} was on the
Intrinsic Low-Dimensional Manifold (ILDM)
method due to Maas and Pope~\cite{MP-1992} and
an iterative method proposed by Fraser~\cite{F-1988}
and further developed by Roussel and Fraser~\cite{RF-1990}.
In this article, the focus is on 
the Computational Singular
Perturbation (CSP) method developed by
Lam and Goussis~\cite{GL-1992, HG-1999,
L-1993, LG-1988, LG-1991, LG-1994,
LJL-2001, MDMG-1999, MG-2001, VG-2001}.

A chemical kinetic equation is an
ordinary differential equation (ODE),
\be
  \frac{dx}{dt} = g(x) ,
  \Label{system}
\ee
for a vector $x$
of species concentrations;
$g$ is a smooth vector field,
and $t$ is time.
Reduction methods are effective
when the variables fall
into two classes, fast and slow,
as is the case when the Jacobian
of the vector field has
a spectral gap.
For the analysis, it is convenient
to identify the spectral gap with
the inverse of a small parameter $\eps$,
but we emphasize that this restriction
is not necessary for the applicability 
of the CSP method.
The characteristic time scales for
the fast and slow species are given by
the ``fast'' time $t$ and
the ``slow'' time $\tau = \eps t$,
respectively.
We assume that the entries of~$x$
are ordered in such a way that
the first $m$ components evolve
on the slow time scale and the remaining
$n$ components on the fast time scale.
Then the vector field $g$ has the form
\be
  g
  =
  \left(
  \begin{array}{c}
    \eps g_1 \\
    g_2
  \end{array}
  \right)
  =
  \left(
  \begin{array}{cc}
    I_m & 0 \\
    0   & I_n
  \end{array}
  \right)
  \left(
  \begin{array}{c}
    \eps g_1 \\
    g_2
  \end{array}
  \right) ,
  \Label{g-12}
\ee
where $I_m$ and $I_n$ are
the identity matrices in
$\mathbf{R}^m$ and $\mathbf{R}^n$, respectively,
and the system~(\ref{system})
is a \emph{fast--slow system} of ODEs.
Both $g_1$ and $g_2$ may depend on $\eps$,
but the entries of these vectors
as well as their partial derivatives
are all $\mathcal{O} (1)$
as $\eps \downarrow 0$,
uniformly in $x$.

Geometric singular perturbation theory
(GSPT)~\cite{F-1979, J-1994} provides
a natural framework for the analysis
of fast--slow systems of ODEs.
If such a system has a slow manifold, $\Mo$,
in the limit as $\eps \downarrow 0$
and this manifold is normally hyperbolic,
then GSPT identifies a (usually nonunique)
slow manifold $\Me$
for $\eps$ sufficiently small.
GSPT also gives a complete geometric and
analytical description of
all solutions near~$\Me$,
including how trajectories
approach~$\Me$.
The goal of any reduction method is
to find $\Me$, if it exists.

Typically, the vector field $g$ is written
in a form suggested by chemical kinetics,
namely, as a weighted sum of the stoichiometric vectors,
the weights being the associated reaction rates.
But this representation is in no way unique.
In fact, Eq.~(\ref{g-12}) shows an equivalent
representation of $g$ as a weighted sum
of the standard basis vectors of $\mathbf{R}^{m+n}$,
the weights being the coordinates
$\eps g_1, \ldots\,, \eps g_m, g_{m+1}, \ldots\,, g_{m+n}$.
The objective of the CSP method is
to express $g$ in yet another basis,
one that is tuned to the dynamics of the system,
where the fast and slow coordinates (amplitudes)
evolve independently of each other.
The CSP method achieves this objective
constructively by successive approximation.
Starting with a more or less arbitrary initial basis,
one derives the evolution equations for
the fast and slow amplitudes
and updates the basis iteratively
in such a way that
the evolution equations for the
updated fast and slow amplitudes
decouple to increasingly higher order
in the small parameter~$\eps$.
Each iteration consists of two steps.
The first step deals with the dependence
of the fast amplitudes on the slow amplitudes,
the second step with the dependence
of the slow amplitudes on the fast amplitudes.

After each iteration,
one identifies the \emph{CSP manifold} (CSPM)
as the locus of points where
the then-current fast amplitudes vanish.
The CSPM is an approximation to the slow manifold $\Me$.
The question is: How good is the approximation?
In this paper, we 
analyze the general class of fast-slow systems
of ODEs~(\ref{system})--(\ref{g-12})
and show (Theorem~\ref{t-CSPMq})
that the CSP method generates term by term
the asymptotic expansion 
of the slow manifold~$\Me$.
After $q$ iterations
($q=0,1,2,\ldots$),
the asymptotic expansions of the
CSPM and~$\Me$ agree up to and including
terms of $\mathcal{O} (\eps^q)$;
they differ in general at $\mathcal{O} (\eps^{q+1})$.
Also, the $q$th application of the CSP algorithm
leaves the terms at $\mathcal{O} (1)$
through $\mathcal{O} (\eps^{q-1})$
invariant.
(This observation is important
because the lower-order terms
have already been determined correctly
in the preceding applications.)
We illustrate Theorem~\ref{t-CSPMq}
with an example from the
Michaelis-Menten-Henri
mechanism of enzyme kinetics~\cite{F-1975,
HTA-1967, O-1991, P-1987, PL-1984}.
Similar results
(for $q = 1, 2$) have been obtained by
Valorani, Goussis, and Najm~\cite{VGN-2002}
for a model equation due to
Davis and Skodje~\cite{DS-1999}.

Our proof proceeds
via an intermediate result
for a one-step CSP method.
The one-step CSP method
is the same as the full
two-step CSP method
but involves only the first step.
It yields a sequence of slow manifolds,
just like the full CSP method,
whose asymptotic behavior as $\eps \downarrow 0$
can be compared with that of the slow manifold~$\Me$.
The result (Theorem~\ref{1step-Main-thm}) is
that $q$ applications of the one-step CSP algorithm
yield an approximate slow manifold that agrees
asymptotically with $\Me$ up to
and including terms of $\mathcal{O} (\eps^q)$.
In other words, the one-step CSP method
is as accurate as the full CSP method;
and to prove the main result
for the full CSP method,
one needs only to show that
the second step does not affect
the lower-order terms in the asymptotic
expansion of the CSPM.
Although the second step of the CSP method
does not play a role in the approximation
of~$\Me$, it does play a constructive role
in the approximation of
the dynamics near $\Me$,
as we shall demonstrate
in the special case of
the Michaelis--Menten--Henri equations.

In~\cite{KK-2002}, we showed
that the ILDM method yields
an approximate slow manifold
that is asymptotically accurate
up to and including terms of
$\mathcal{O} (\eps)$,
with an error of $\mathcal{O} (\eps^2)$
proportional to the curvature of~$\Mo$.
The CSP method, on the other hand,
can generate an approximate slow manifold
that is asymptotically accurate
up to any order.
The difference can be traced to two facts,
namely,
the choice of the fundamental operator
governing the dynamics of the system
and the retention of the variation of
the Jacobian over the manifold~$\Mo$.
While the ILDM method is designed to transform
the Jacobian of the vector field
into triangular form 
(and often also into diagonal form),
the CSP method is an algorithm
to diagonalize the (nonlinear)
Lie bracket involving the vector field
to successively higher orders in $\eps$.
The Jacobian is a linear approximation,
so the ILDM method never gets beyond
a linear approximation.
The variation of the Jacobian
over $\Mo$ introduces an extra term
in the Lie bracket.
By retaining it, the CSP method
preserves the nonlinear character
of the operator governing the dynamics
of the system.

This article
is organized as follows.
In Section~\ref{s-general},
we recall the Fenichel theory
of GSPT and give the asymptotic expansion
of the slow manifold~$\Me$.
In Section~\ref{s-CSP},
we describe the full CSP method
for fast--slow systems and
state Theorem~\ref{t-CSPMq}.
The one-step CSP method is introduced
in Section~\ref{s-1step}.
The approximation result for
the slow manifold is given
in Theorem~\ref{1step-Main-thm};
its proof occupies most of Section~\ref{s-1step}
and uses two lemmas that are given
in the Appendix.
In Section~\ref{s-Full},
we return to the full CSP method
and prove Theorem~\ref{t-CSPMq}.
In Section~\ref{s-MMH},
we illustrate the CSP method and the results
of this paper on a planar system of equations
for the Michaelis--Menten--Henri mechanism
of enzyme kinetics.
Section~\ref{s-discussion}
is devoted to a discussion of
the relation between the CSP
and ILDM methods.

\section{Fast-Slow Systems of ODEs \label{s-general}}
\setcounter{equation}{0}
Collecting the slow variables
in a single (column) vector $y$
and the fast variables
in a (column) vector $z$,
we rewrite Eqs.~(\ref{system})--(\ref{g-12})
as a fast-slow system,
\begin{eqnarray}
  y' &=& \eps g_1 (y, z, \eps) ,  \Label{eq-y} \\
  z' &=& g_2 (y, z, \eps) . \Label{eq-z}
\end{eqnarray}
(A prime ${}'$ denotes differentiation
with respect to $t$.)
The long-term dynamics of this system
are more naturally studied on the time scale
of the slow variable $\tau = \eps t$,
where the system of Eqs.~(\ref{eq-y})--(\ref{eq-z})
assumes the form
\begin{eqnarray}
  \dot{y} &=& g_1 (y, z, \eps) ,  \Label{eq-y-slow} \\
  \eps \dot{z} &=& g_2 (y, z, \eps) . \Label{eq-z-slow}
\end{eqnarray}
(A dot $\dot{\ }$ denotes differentiation
with respect to $\tau$.)

In the limit $\eps \downarrow 0$,
Eq.~(\ref{eq-z-slow}) reduces formally
to the algebraic equation
$g_2 (y, z, 0) = 0$.
We assume that there exists
a compact domain $K \in \mathbf{R}^m$
and a smooth single-valued
function $h_0$ on $K$
such that
\be
  g_2 (y, h_0(y), 0) = 0 , \quad y \in K .
  \Label{h0}
\ee
Then the long-time dynamics
of the system~(\ref{eq-y})--(\ref{eq-z})
are confined to the 
\emph{reduced slow manifold}~$\Mo$,
\be
  \mathcal{M}_0 = \{ (y, z) \in \mathbf{R}^{m+n} :
  z = h_0 (y), \; y \in K \} .
  \Label{M0}
\ee
We assume, furthermore,
that the real parts of the eigenvalues
of the matrix $D_z g_2 (y, h_0(y), 0)$
are all negative,
so $\Mo$ is asymptotically stable.
Then the Fenichel theory~\cite{F-1979},
which applies more generally to
normally hyperbolic invariant manifolds,
guarantees that $\Mo$ persists as
a slow manifold, so
for all sufficiently small $\eps$
there exists a \emph{slow manifold}, $\Me$ ,
that is invariant under the dynamics
of the system of Eqs.~(\ref{eq-y})--(\ref{eq-z}).
Moreover, $\Me$ has the same dimension as $\Mo$
and lies near $\Mo$,
all nearby solutions relax
exponentially fast to $\Me$,
and the long-term dynamics
of the system~(\ref{eq-y})--(\ref{eq-z})
are governed by an equation
on $\Me$.
The manifold $\Me$
is not unique;
typically, there is a family
of slow manifolds, 
all exponentially close
($\mathcal{O} (\mathrm{e}^{-c/\eps})$ 
for some $c>0$).
The following theorem is essentially a
restatement of~\cite[Theorem 2]{J-1994}.

\begin{THEOREM} \Label{t-Fenichel}
For all sufficiently small $\eps$,
there is a function $h_\eps$
such that the graph
\be
  \Me = \{ (y,z) : z = h_\eps (y) , \; y \in K \}
  \Label{M-eps}
\ee
is locally invariant 
under the dynamics of
Eqs.~(\ref{eq-y})--(\ref{eq-z}).
The function $h_\eps$ admits
an asymptotic expansion
as $\eps \downarrow 0$,
\be
  h_\eps (y)
  = h_0 (y) + \eps h_1 (y)
  + \eps^2 h_2 (y) + \cdots ,
  \Label{h-eps-exp}
\ee
and $h_\eps \in C^r (K)$
for any finite $r$.
The long-term dynamics 
of the system  of
Eqs.~(\ref{eq-y})--(\ref{eq-z})
are governed by  the equation
\be
  \dot{y} = g_1 ( y, h_\eps(y), \eps)
  \Label{traj}
\ee
on $\Me$,
where $\dot{\ } = d/d\tau$ with $\tau = \eps t$.
\end{THEOREM} 

The coefficients $h_1, h_2, \ldots$
are found from the \emph{invariance equation},
\be
  g_2 (y, h_\eps (y), \eps)
  - \eps (Dh_\eps) (y) g_1 (y, h_\eps (y), \eps)
  = 0 , \quad y \in K ,
  \Label{inveq}
\ee
in the following manner.
(The invariance equation follows immediately
from the chain rule, $z'=Dh_\eps(y) y'$,
and Eqs.~(\ref{eq-y})--(\ref{eq-z}).)
Each of the functions
$g_1 ( \cdot\,, h_\eps, \eps)$
and
$g_2 ( \cdot\,, h_\eps, \eps)$
admits a Taylor expansion near $\eps=0$,
\be
  g_1 (\cdot\,, h_\eps, \eps)
  =
  \sum_{q=0}^\infty g_{1,q} \eps^q , \quad
  g_2 (\cdot\,, h_\eps, \eps)
  =
  \sum_{q=0}^\infty g_{2,q} \eps^q ,
\ee
with coefficients
\begin{eqnarray}
    g_{1,q}
    &=&
    \sum_{k=0}^{q-1} \sum_{j=1}^{q-k} \frac{1}{k!j!} (D^j_z D^k_\eps g_1)_0 
    \sum_{\vert i\vert=q-k} (h_{i_1},...,h_{i_j})
    + \frac{1}{q!} (D^q_\eps g_1)_0 , \Label{Ho2oH1} \\
    g_{2,q}
    &=&
    \sum_{k=0}^{q-1} \sum_{j=1}^{q-k} \frac{1}{k!j!} (D^j_z D^k_\eps g_2)_0 
    \sum_{\vert i\vert=q-k} (h_{i_1},...,h_{i_j})
    + \frac{1}{q!} (D^q_\eps g_2)_0 . \Label{Ho2oH2}
\end{eqnarray}
The notation $(\,\cdot\,)_0$
indicates that the quantity inside the parentheses
is evaluated on $\Mo$---that is,
at $(y, h_0(y), 0)$.
Note that
$(D_z^j D_\eps^k g)$
is a multilinear operator,
which maps a $j$-form to a vector.
The inner sum in Eqs.~(\ref{Ho2oH1}) and
(\ref{Ho2oH2}) is taken over
all multiindices $i = (i_1, \ldots\,, i_j)$
of $j$ positive integers $i_1$ through $i_j$
subject to the constraint
$|i| = i_1 + \cdots\, + i_j = q-k$.
The expressions~(\ref{Ho2oH1}) and (\ref{Ho2oH2})
hold for all $q$ if it is understood
that a sum is empty whenever
its lower bound exceeds its upper bound.
Substituting the expansions~(\ref{Ho2oH1})
and (\ref{Ho2oH2}) into the invariance equation,
Eq.~(\ref{inveq}),
and setting the coefficient of $\eps^q$ equal to zero,
we obtain an infinite set of equations,
\be
  g_{2,q}
  - \sum_{\ell=0}^{q-1} (Dh_\ell) g_{1,q-1-\ell}
  = 0 , \quad q = 0, 1, \ldots\,.
  \Label{inveq-q}
\ee
The first few equations are
\begin{eqnarray}
  &&g_{2,0} = 0 , \Label{h(0)} \\
  &&(D_z g_2)_0 h_1
  + (D_\eps g_2)_0
  - (Dh_0) g_{1,0}  = 0 , \Label{h(1)} \\
  &&(D_z g_2)_0 h_2
  + \half (D_z^2 g_2)_0 \left( h_1, h_1 \right)
  + (D_z D_\eps g_2)_0 h_1
  + \half (D_\eps^2 g_2)_0 \nonumber \\
  &&\hspace{1em}\mbox{}
  - (D h_1) g_{1,0}
  - (Dh_0) \left( (D_z g_1)_0 h_1
  - (D_\eps g_1)_0 \right) = 0 . \Label{h(2)}
\end{eqnarray}
Equation~(\ref{h(0)}) is satisfied identically,
Eq.~(\ref{h(1)}) yields the coefficient $h_1$,
Eq.~(\ref{h(2)}) the coefficient $h_2$, and so on.

\vspace{-2ex}\paragraph{Remark.}
The assumption that the chemical species
can be divided into fast and slow species,
as in Eqs.~(\ref{eq-y})--(\ref{eq-z}),
is made for convenience.
Our analysis can also be applied 
to general chemical systems
where each species may be involved
in both fast and slow reactions
and for which there is a slow manifold.

\section{The CSP Method for Fast--Slow Systems \label{s-CSP}}
\setcounter{equation}{0}
In Eq.~(\ref{g-12}), the vector field $g$
is represented in terms of 
the standard orthonormal basis.
It is useful to examine the representation of
$g$ in terms of other bases,
especially bases whose entries depend on $x$.

Let $A$ be an $(m+n) \times (m+n)$ matrix
whose entries may depend on $x$
and
whose columns form a basis for
the space $\mathbf{R}^{m+n}$
for each $x$.
The vector field $g$ 
may be expressed in terms of 
this (variable) basis $A$ as
\be
  g = A f ,
\Label{g=Af}
\ee
where $f$ is the vector of
the coordinates (amplitudes) of $g$.
(When the columns of $A$
are the stoichiometric vectors,
the amplitudes are
the reaction rates.)
Since $A$ is invertible,
$B = A^{-1}$, and
\be
  f = B g .
\Label{f}
\ee
The amplitudes can be split
into two classes,
$f =
\left( \begin{array}{c} 
      f^1 \\ f^2 
       \end{array}
\right)$,
where $f^1$ is an $n$-vector 
representing the fast amplitudes
and $f^2$ an $m$-vector
representing the slow amplitudes.
The splitting suggests
that we consider a decomposition
of $A$, namely,
$A =
  \left( A_1, A_2 \right)$,
where 
$A_1$ is $(m+n) \times n$
and $A_2$ is $(m+n) \times m$,
and a corresponding decomposition
of $B$, namely,
$B =
\left( \begin{array}{c} 
       B^1 \\ B^2 
       \end{array}
\right)$,
where $B^1$ is $n \times (m+n)$,
and $B^2$ is $m \times (m+n)$.
Thus,
$f^1 = B^1 g$ and $f^2 = B^2 g$.
Also, the identity $AB=I$
on $\mathbf{R}^{m+n}$
implies that
$A_1 B^1 + A_2 B^2 = I$,
while the identity $BA=I$
on $\mathbf{R}^{m+n}$
implies that
$B^1 A_1 = I_n$ and
$B^2 A_1 = 0$ on $\mathbf{R}^n$
and
$B^2 A_2 = I_m$ and
$B^1 A_2 = 0$ on $\mathbf{R}^m$.

The fast and slow amplitudes
evolve in time.
Differentiating Eq.~(\ref{f})
along solutions of the system~(\ref{system}), 
we obtain
\[
  \frac{df}{dt}
  =
  B \frac{dg}{dt} + \frac{dB}{dt} g
  =
  B (Dg) g + \frac{dB}{dt} g ,
\]
where $Dg$ is the Jacobian of $g$.
Hence, $f$ satisfies the nonlinear ODE
\be
  \frac{df}{dt}
  = 
  \Lambda f,
\ee
where $\Lambda$, the generator of the dynamics
for the amplitudes, is given by
\be
  \Lambda
  =
  B (Dg) A + \frac{dB}{dt} A .
  \Label{Lambda}
\ee
Since $BA = I$ and $I$ is time invariant,
$A$, $B$, and their time derivatives
satisfy the identity
\be
  (dB/dt)A + B(dA/dt) = 0
  \Label{BA-inv}
\ee
at all times.
Hence, the definition~(\ref{Lambda})
is equivalent to
\be
  \Lambda
  =
  B (Dg) A - B \frac{dA}{dt} ,
  \Label{Lambda-alt}
\ee
where $dA/dt = (DA)g$.
For completeness,
we note that
the identity~(\ref{BA-inv})
implies that
$((DB)Af)A + B((DA)Af) = 0$.

In general, the operator $\Lambda$ is not diagonal,
and the equations governing the evolution of
$f^1$ and $f^2$ are coupled. 
An ideal basis $A$
is one in which $\Lambda$
is block-diagonalized,
so that the ODEs for
$f^1$ and $f^2$ decouple.
The CSP method approaches this ideal
by successive refinements of
the basis matrices $A$ and $B$.
The algorithm starts from
a constant matrix $\Af^{(0)}$,
\be
  \Af^{(0)}
  =
  \left( \Af_1^{(0)}, \, \Af_2^{(0)} \right)
  =
  \left( 
  \begin{array}{cc}
    \Af_{11}^{(0)} & \Af_{12}^{(0)} \\
    \Af_{21}^{(0)} & \Af_{22}^{(0)}
  \end{array}
  \right) .
  \Label{A(0)-full}
\ee
Here,
$\Af_{11}^{(0)}$ is an $m \times n$ matrix,
$\Af_{22}^{(0)}$ an $n \times m$ matrix,
and the off-diagonal blocks
$\Af_{12}^{(0)}$ and $\Af_{21}^{(0)}$
are full-rank square matrices of order
$m$ and $n$, respectively.
A common choice is
$\Af_{11}^{(0)} = 0$,
so every column vector of $\Af_1^{(0)}$ lies
in the fast subspace.
We follow this convention and assume,
henceforth, that $\Af_{11}^{(0)} = 0$,
\be
  \Af^{(0)}
  =
  \left( \Af_1^{(0)}, \, \Af_2^{(0)} \right)
  =
  \left( 
  \begin{array}{cc}
    0            & \Af_{12}^{(0)} \\
    \Af_{21}^{(0)} & \Af_{22}^{(0)}
  \end{array}
  \right) .
  \Label{A(0)}
\ee
A more general choice of $A^{(0)}$
is discussed below, after Theorem~\ref{t-CSPMq}.
The inverse of $\Af^{(0)}$ is
\begin{eqnarray}
  \Bf_{(0)}
  &=& ( \Af^{(0)} )^{-1}
  = \left(
  \begin{array}{c} \Bf^1_{(0)} \\ \Bf^2_{(0)} \end{array}
  \right)
  =
  \left( 
  \begin{array}{cc}
    \Bf^{11}_{(0)} & \Bf^{12}_{(0)} \\
    \Bf^{21}_{(0)} & 0
  \end{array}
  \right) \nonumber \\
  &=& 
  \left(
  \begin{array}{cc}
    - ( \Af_{21}^{(0)} )^{-1}
    \Af_{22}^{(0)}
    ( \Af_{12}^{(0)} )^{-1}
    & ( \Af_{21}^{(0)} )^{-1} \\
    ( \Af_{12}^{(0)} )^{-1}
    & 0
  \end{array}
  \right) .
  \Label{B(0)}
\end{eqnarray}

The algorithm then proceeds iteratively.
For $q = 0, 1, \ldots\,$,
one first defines the matrix
$\bfLf_{(q)}$
in accordance with
Eq.~(\ref{Lambda-alt}),
\be
  \bfLf_{(q)}
  =
  \Bf_{(q)} (Dg) \Af^{(q)} - \Bf_{(q)} \frac{d\Af^{(q)}}{dt}
  =
  \left(
  \begin{array}{cc}
    \bfLf_{(q)}^{11} & \bfLf_{(q)}^{12} \\
    \bfLf_{(q)}^{21} & \bfLf_{(q)}^{22}
  \end{array}
  \right) ,
  \Label{Lambda(q)}
\ee
and matrices
$\Uf_{(q)}$
and
$\Lf_{(q)}$,
\be
  \Uf_{(q)}
  =
  \left(
  \begin{array}{cc}
    0 & (\bfLf_{(q)}^{11})^{-1} \bfLf_{(q)}^{12} \\
    0 & 0
  \end{array}
  \right) , \quad
  \Lf_{(q)}
  =
  \left(
  \begin{array}{cc}
    0 & 0 \\
    \bfLf_{(q)}^{21} (\bfLf_{(q)}^{11})^{-1} & 0
  \end{array}
  \right) .
\Label{U(q),L(q)}
\ee
Then one updates
$A^{(q)}$ and $B_{(q)}$
according to the formulas
\begin{eqnarray}
  \Af^{(q+1)}
  &=&
  \Af^{(q)} ( I - \Uf_{(q)} ) ( I + \Lf_{(q)} ) , \Label{A(q+1)} \\
  \Bf_{(q+1)}
  &=&
  ( I - \Lf_{(q)} ) ( I + \Uf_{(q)} ) \Bf_{(q)} ,  \Label{B(q+1)}
\end{eqnarray}
and returns to Eq.~(\ref{Lambda(q)})
for the next iteration.

At each iteration, one imposes
the \emph{CSP condition},
\be
  B_{(q)}^1 g = 0 , \quad q = 0, 1, \ldots \,,
\Label{CSPM-q}
\ee
to identify those points
where the fast reaction rates
vanish with respect to
the then-current basis.
For $q=0$, $B_{(0)}^1$ is constant
and given by Eq.~(\ref{B(0)});
for $q = 1, 2, \ldots\,$,
the CSP condition takes the form
\be
  B^1_{(q)}(y,\psi_{(q-1)}(y,\eps),\eps)
  g(y,z,\eps) = 0 , \quad
  q = 1, 2, \ldots \,.
\ee
If, for any $q$,
the CSP condition
is satisfied by a function
$z = \psif_{(q)}(y,\eps)$,
then
\be
  \Kf^{(q)}
  = \{ (y,z): z = \psif_{(q)} (y, \eps) , \; y \in K \} ,
  \quad q = 0, 1, \ldots
\Label{L-q}
\ee
is defined as the
\emph{CSP manifold} (CSPM) of order $q$.

\begin{THEOREM} \Label{t-CSPMq}
The CSP manifold $\Kf^{(q)}$
agrees asymptotically with~$\Me$
up to and including terms
of $\mathcal{O} (\eps^q)$
for $q = 0, 1, \ldots\,$,
\be
\psif_{(q)} ( \cdot\,, \eps)
=
\sum_{j=0}^q \eps^j h_j + \mathcal{O} (\eps^{q+1}) , \quad
\eps \downarrow 0 .
\ee
\end{THEOREM}

Our proof of Theorem~\ref{t-CSPMq}
proceeds via an intermediate result,
which is of independent interest.
We introduce a ``truncated'' CSP method,
where we apply, at each iteration,
only the first of the two steps
of the full CSP method
and skip the second step.
This \emph{one-step CSP method}
reduces the matrix $\Lambda$
to lower block-triangular form.
We show that,
after $q$ iterations,
the one-step CSP method generates
a manifold $\Ks^{(q)}$,
whose asymptotic expansion
agrees with that of $\Me$
up to and including terms of
$\mathcal{O}(\eps^q)$
(Theorem~\ref{1step-Main-thm}).
In other words, the one-step
CSP method is as accurate as
the full CSP method is claimed
to be in Theorem~\ref{t-CSPMq}.
We then return to the full CSP method
and carry out an asymptotic analysis
of the modifications introduced by
the second step.
This second step reduces $\Lambda$
further to block-diagonal form.
We show that,
at the $q$th iteration,
the second step affects only terms of
$\mathcal{O}(\eps^{q+1})$ and higher.
Hence, $\Kf^{(q)}$ approximates
$\Me$ as accurately as
$\Ks^{(q)}$, and
Theorem~\ref{t-CSPMq} follows.

Theorem~\ref{t-CSPMq}
extends readily to the case
where the eigenvectors of the Jacobian~$Dg$
are used, instead of the stoichiometric vectors,
to form the initial basis~$A^{(0)}$.
In that case,
the slow subspace of
the leading-order Jacobian
coincides with the tangent space
${\mathcal T}_p \Mo$ at any point
$p \in \Mo$,
so the columns of $A_2^{(0)}$
are tangent to $\Mo$
to leading order.
In turn, this implies
that the rows of $\Bf^1_{(0)}(p)$ span
the orthogonal complement of the tangent space,
also to leading order.
As a result, the initial CSPM,
the solution of $B^1_{(0)}g=0$,
coincides with $\Me$ 
up to and including terms
of $\mathcal{O}(\eps)$,
which is one order
higher than is the case
when $A^{(0)}$ is given by Eq.~(\ref{A(0)}).
Moreover, for each $q=1,2,\ldots\,$,
the proof of Theorem~\ref{1step-Main-thm}
generalizes directly to this case. 
The asymptotic expansion of $\psis_{(q)}$ 
coincides with that of $h_\eps$
up to and including terms of $\mathcal{O}(\eps^{q+1})$,
which is one order
higher than is the case
when $A^{(0)}$ is given
by Eq.~(\ref{A(0)}).

\paragraph{Remark.}
Lam and Goussis,
in their presentation~\cite{L-1993}
of the CSP method,
perform the update~(\ref{A(q+1)})
and~(\ref{B(q+1)}) in two steps.
The first step corresponds to the postmultiplication
of $\Af^{(q)}$ with $I - \Uf_{(q)}$ and
premultiplication of $\Bf_{(q)}$ with $I + \Uf_{(q)}$,
the second step to the subsequent postmultiplication
of $\Af^{(q)}(I - \Uf_{(q)})$ with $I + \Lf_{(q)}$ and
premultiplication of $(I + \Uf_{(q)})\Bf_{(q)}$
with $I - \Lf_{(q)}$.
The nonzero entries of $U_{(q)}$ and $L_{(q)}$
are chosen so that $\Lambda$
is block-diagonalized
to successively
higher order in $\eps$.

\paragraph{Remark.}
The definition~(\ref{Lambda-alt})
implies that $\Lambda$ is
the product of $B$
with the Lie bracket of $A$
(considered column by column)
and $g$,
\be
  \Lambda
  = B \left[A,g\right]
  = B ([\Af_{\,\cdot\,,1}, g], \ldots, [\Af_{\,\cdot\,,m+n}, g]) .
\Label{L=Lie bracket}
\ee
The Lie bracket of two vector fields $a$ and $g$
is $[a, g] = (Dg) a - (Da) g$~\cite{O-1986}.

\paragraph{Remark.}
It is useful to
state how
$\Lambda$ transforms
to understand its properties as
an operator.
If $\hat{A} = AC$ and
$\hat{B} = C^{-1} B$,
where 
$C$ is an invertible square matrix
representing a coordinate transformation
in $\mathbf{R}^{m+n}$,
then
\begin{eqnarray}
  \hat{\Lambda}
  &=&
  \hat{B} (Dg) \hat{A} - \hat{B} \frac{d\hat{A}}{dt}
  =
  C^{-1} B (Dg) AC - C^{-1} B \frac{d(AC)}{dt} \nonumber \\
  &=&
  C^{-1} B (Dg) AC
  - C^{-1} B \left( \frac{dA}{dt} C + A \frac{dC}{dt} \right) \nonumber \\
  &=&
  C^{-1} \Lambda C - C^{-1} \frac{dC}{dt} ,
\Label{L-transf}
\end{eqnarray}
where $dC/dt = (DC)g$.
The presence of
the term
$C^{-1} dC/dt$
in Eq.~(\ref{L-transf})
shows that
$\hat{\Lambda}$
and
$\Lambda$
are not similar
unless $C$ is constant.

\section{The One-Step CSP Method \label{s-1step}}
\setcounter{equation}{0}
The goal of the one-step CSP method is
to reduce the matrix $\Lambda$
to lower block-triangular form---that is,
to push the matrix $\Lambda^{12}$
to increasingly higher order in $\eps$.
The method is identical to
the full CSP method
except for the updating of the matrices $A$ and $B$.
One starts from
the same bases, 
$\As^{(0)}=\Af^{(0)}$
and
$\Bs_{(0)}=\Bf_{(0)}$,
and,
instead of Eqs.~(\ref{A(q+1)}) and~(\ref{B(q+1)}),
uses the one-step expressions
\begin{eqnarray}
    \As^{(q+1)} &=& \As^{(q)} (I-\Us_{(q)}) , \Label{A-rec} \\
    \Bs_{(q+1)} &=& (I+\Us_{(q)}) \Bs_{(q)} , \Label{B-rec}
\end{eqnarray}
where the matrix $\Us_{(q)}$
is defined as in Eq.~(\ref{U(q),L(q)})
with $\bfLf$ replaced by $\bfLs$.
(A tilde $\tilde{\ }$
distinguishes
a quantity
from its counterpart
in the full CSP method.)

The update rule for $\bfLs$
follows immediately from Eq.~(\ref{L-transf}),
\be
  \bfLs_{(q+1)}
  =
  (I+\Us_{(q)}) \bfLs_{(q)} (I-\Us_{(q)})
  + (I+\Us_{(q)}) \frac{d\Us_{(q)}}{dt} .
\Label{L-rec-old}
\ee
(Note that the identities
$\As^{(0)}=\Af^{(0)}$
and
$\Bs_{(0)}=\Bf_{(0)}$
imply that
$\bfLs_{(0)} = \bfLf_{(0)}$.)
The matrix $\Us_{(q)}$
and its time derivative
have the same block structure;
only the upper right block
is nonzero,
so $\Us_{(q)} d\Us_{(q)}/dt = 0$,
and Eq.~(\ref{L-rec-old}) reduces to
\be
  \bfLs_{(q+1)}
  =
  (I+\Us_{(q)}) \bfLs_{(q)} (I-\Us_{(q)})
  + \frac{d\Us_{(q)}}{dt} .
\Label{L-rec}
\ee
In terms of the constituent blocks, we have
\begin{eqnarray}
    \bfLs^{11}_{(q+1)} &=& \bfLs^{11}_{(q)} + \Us_{(q)} \bfLs^{21}_{(q)}, 
\Label{L11-rec}\\
    \bfLs^{12}_{(q+1)} &=& \Us_{(q)} \bfLs^{22}_{(q)}
                          - \Us_{(q)} \bfLs^{21}_{(q)} \Us_{(q)} + \frac{d\Us_{(q)}}{dt} ,
\Label{L12-rec}\\
    \bfLs^{21}_{(q+1)} &=& \bfLs^{21}_{(q)},
\Label{L21-rec}\\
    \bfLs^{22}_{(q+1)} &=& \bfLs^{22}_{(q)} - \bfLs^{21}_{(q)} \Us_{(q)}, 
\Label{L22-rec}
\end{eqnarray}
where we have used Eq.~(\ref{U(q),L(q)})
to simplify Eq.~(\ref{L12-rec}).
Note that we freely use $\Us_{(q)}$ to denote
both the full update matrix
and its restriction to the subspace $\mathbf{R}^m$;
the latter is represented by the matrix
$(\bfLs_{(q)}^{11})^{-1} \bfLs_{(q)}^{12}$.
The appropriate interpretation is clear from the context.

The one-step CSP method
generates a sequence of manifolds,
\be
  \Ks^{(q)}
  = \{ (y,z): z = \psis_{(q)} (y, \eps) , \; y \in K \} ,
  \quad q = 0, 1, \ldots\, ,
\Label{K-q}
\ee
just like the full CSP method;
cf.~Eq.~(\ref{L-q}).
The functions $\psis_{(q)}$
are defined by the conditions
\be
  \Bs_{(q)}^1 g = 0 , \quad q = 0, 1, \ldots\, ,
  \Label{CSPM-q-1step}
\ee
where $\Bs_{(q)}^1$
is obtained from Eq.~(\ref{B-rec}).

\begin{THEOREM} \Label{1step-Main-thm}
The manifold $\Ks^{(q)}$
agrees asymptotically with~$\Me$
up to and including terms
of $\mathcal{O} (\eps^q)$
for $q = 0, 1, \ldots\,$,
\be
\psis_{(q)} ( \cdot\,, \eps)
=
\sum_{j=0}^q \eps^j h_j + \mathcal{O} (\eps^{q+1}) , \quad
\eps \downarrow 0 .
\ee
\end{THEOREM}
The proof of the theorem is by induction on $q$.

\subsection{The Induction Hypothesis
\label{ss-strategy}}
The central idea of the proof of
Theorem~\ref{1step-Main-thm}
is to express
the CSP condition~(\ref{CSPM-q-1step})
in a form that resembles
that of the invariance equation~(\ref{inveq})
and then to derive the conditions
under which
the left and right members of
the two equations
are the same at each order.

We begin by expressing
the quantities
$\As^{(q+1)}$,
$\Bs_{(q+1)}$,
and
$\bfLs_{(q+1)}$
in terms of
the original quantities
$\Af^{(0)}$,
$\Bf_{(0)}$,
and
$\bfLf_{(0)}$.
Applying the definition~(\ref{A-rec})
recursively, we find
\[
\As^{(q+1)}
=
\Af^{(0)} \prod_{j=0}^q (I - \Us_{(j)}) .
\]
Since each $\Us_{(j)}$ is nilpotent,
it follows that
\be
  \As^{(q+1)} = \Af^{(0)} (I - \Ps_{(q)}) ,
\Label{A-expl}
\ee
where
\be
  \Ps_{(q)}
  = \sum_{j=0}^q \Us_{(j)}
  = \left(
  \begin{array}{cc}
  0 & \sum_{\ell=0}^q 
  (\bfLs^{11}_{(\ell)})^{-1} \bfLs^{12}_{(\ell)} \\
  0 & 0
  \end{array}
  \right) .
\Label{expl-aux1}
\ee
Similarly,
\be
  \Bs_{(q+1)} = (I + \Ps_{(q)}) \Bf_{(0)} .
\Label{B-expl}
\ee
Substituting Eqs.~(\ref{A-expl}) and (\ref{B-expl})
into the transformation formula~(\ref{L-transf}),
and recalling that 
$\bfLs_{(0)} = \bfLf_{(0)}$
and $\Ps_{(q)}d\Ps_{(q)}/dt=0$, 
we find
\be
  \bfLs_{(q+1)} = (I + \Ps_{(q)}) \bfLf_{(0)} (I - \Ps_{(q)}) 
             + \frac{d\Ps_{(q)}}{dt} .
\Label{L-expl}
\ee
We use these expressions to rewrite
Eq.~(\ref{CSPM-q-1step}).
Since $\Bf_{(0)}^{22}=0$,
the equation becomes
\[
  \Bf_{(0)}^{12} g_2
  +
  \eps
  \left[\Ps_{(q-1)} \Bf_{(0)}^{21} + \Bf_{(0)}^{11}\right] g_1
  = 0
\]
or, since $\Bf^{12}_{(0)} = (\Af^{(0)}_{21})^{-1}$,
\be
  g_2
  +
  \eps \Af^{(0)}_{21}
  \left[
  \Ps_{(q-1)} \Bf_{(0)}^{21} + \Bf_{(0)}^{11}
  \right] g_1
  = 0 .
\Label{CSP-cond(q)}
\ee
The last equation has the same form
as the invariance equation~(\ref{inveq}).
The solution of Eq.~(\ref{inveq}) is
$z = h_\eps (y)$,
which defines $\Me$,
while the solution of Eq.~(\ref{CSP-cond(q)})
is $z = \psis_{(q)} (y, \eps)$,
which defines $\Ks^{(q)}$.

We analyze
the CSP condition~(\ref{CSP-cond(q)})
order by order,
up to and including
the terms of
$\mathcal{O}(\eps^q)$.
We recall that the components
of the vector field
$g (y,z,\eps)$
are evaluated at
$z=\psis_{(q)}(y,\eps)$,
the matrix $\Ps_{(q-1)}$
is evaluated at
$z=\psis_{(q-1)}(y,\eps)$, and
the blocks of
$\Af^{(0)}$
and
$\Bf_{(0)}$
are constant.
Substituting the asymptotic expansion of
$\psis_{(q)}$,
\be
  \psis_{(q)} (y, \eps)
  = \sum_{j=0}^\infty \eps^j \psis_{(q,j)}(y) , \quad
  \eps \downarrow 0 ,
\Label{psi-exp}
\ee
into Eq.~(\ref{CSP-cond(q)})
and setting the coefficients 
of $1, \eps, \ldots\,, \eps^q$
equal to zero,
we obtain a set of equations,
\be
  g_{2,j}
  + \Af^{(0)}_{21}
  \left[
  \Ps_{(q-1,0)} \Bf_{(0)}^{21} + \Bf_{(0)}^{11}
  \right]
  g_{1,j-1}
  +
  \sum_{\ell=1}^{j-1}
  \Af^{(0)}_{21} \Ps_{(q-1,\ell)} \Bf_{(0)}^{21} g_{1,j-\ell-1}
  = 0 ,
\Label{CSP-cond(q,j)}
\ee
for $j=0,1,\ldots\,,q$.
Here,
$\Ps_{(q-1,\ell)}$
is the coefficient of
the $\mathcal{O}(\eps^\ell)$
term in
the asymptotic expansion of
$\Ps_{(q-1)}$.

Equation~(\ref{CSP-cond(q,j)})
defines $\psis_{(q,j)}$
for $j = 0, 1, \ldots\,, q$.
The leading-order ($j=0$) equation 
in the system~(\ref{CSP-cond(q,j)})
is the same for all $q$,
\be
    g_2 (y, \psis_{(q,0)} (y), 0) = 0 , \quad q=0,1,\ldots\, .
\Label{CSP-cond(q,0)}
\ee
This is also the equation defining $h_0$.
Its solution need not be unique,
but we can identify each
$\psis_{(q,0)}$ with $h_0$,
\be
  \psis_{(q,0)}(y) = h_0(y) , \quad q = 0,1,\ldots\, .
\Label{psi(q,0)}
\ee
Then also $\psis_{(q)} (\cdot\,, 0) = h_0$
for $q=0,1,\ldots\,$,
so to leading order
each manifold $\Ks^{(q)}$
coincides with $\Mo$.

We wish to show that
$\psis_{(q,j)}=h_j$
also for $j=1,2,\ldots\,,q$.
To this end, we compare
Eqs.~(\ref{inveq-q})
and~(\ref{CSP-cond(q,j)}).
For a fixed $j$,
the two equations match if
\begin{eqnarray}
    \Af^{(0)}_{21} 
    \left[\Ps_{(q-1,0)} \Bf_{(0)}^{21} + \Bf_{(0)}^{11}\right] 
    &=& 
    - Dh_0 , \Label{Dho-eq} \\
    \Af^{(0)}_{21} \Ps_{(q-1,\ell)} \Bf_{(0)}^{21} 
    &=& 
    - Dh_\ell,  \quad
    \ell = 1, \ldots\,, j-1 .
\Label{Dhl-eq}
\end{eqnarray}
Conversely,
if Eqs.~(\ref{Dho-eq})
and
(\ref{Dhl-eq})
hold,
then
$\psis_{(q,j)}=h_j$.
Notice that
Eqs.~(\ref{Dho-eq})
and~(\ref{Dhl-eq})
are independent of $j$;
hence, they are nested,
in the sense that,
when $j$ is increased by one,
the equations
for lower values of $j$
remain the same.
Thus, it suffices
to prove
Eqs.~(\ref{Dho-eq})
and~(\ref{Dhl-eq})
for $j=q$.
The proof is by induction on $q$,
where the induction hypothesis is
\begin{eqnarray}
    \Us_{(q-1)}(\cdot\,,\psis_{(q-1)},\eps) 
    &=&
   \mathcal{O}(\eps^{q-1}) ,
\Label{U(q-1)-order}\\
    \Af^{(0)}_{21} 
    \left[\Ps_{(q-1)}(\cdot\,,\psis_{(q-1)},\eps) 
          \Bf_{(0)}^{21} + \Bf_{(0)}^{11}\right] 
    &=& 
    - \sum_{j=0}^{q-1} \eps^j Dh_j 
    + \mathcal{O}(\eps^q) ,
\Label{SumU-assum}\\
   \psis_{(q)} (\cdot\,,\eps)
   &=& 
    \sum_{j=0}^q \eps^j h_j + \mathcal{O} (\eps^{q+1}) .
\Label{psi(q)-exp}
\end{eqnarray}
The validity of these equations for $q=1$
is shown in Section~\ref{ss-1step-q=1}.
The induction step is carried out
in Section~\ref{ss-1step-ind}.

\subsection{Proof of Theorem~\ref{1step-Main-thm} for $\mathbf{q=1}$
\label{ss-1step-q=1}}
We fix $q=1$
and consider the $\mathcal{O}(\eps)$
terms of Eq.~(\ref{CSP-cond(q)}),
\be
      (D_z g_2)_0 \psis_{(1,1)} 
    + (D_\eps g_2)_0 + \Af^{(0)}_{21} \left[
      \Ps_{(0,0)} \Bf_{(0)}^{21} + \Bf_{(0)}^{11}
                                    \right] g_{1,0} = 0.
\Label{CSP-cond(1,1)-aux}
\ee
The first and second terms in this equation
are exactly the same as those in
the equation for $h_1$, see (\ref{h(1)}).
Therefore, we need only
to show that the third term
equals $-(Dh_0)g_{1,0}$
in order to prove the theorem for $q=1$.

According to the definitions~(\ref{expl-aux1})
and~(\ref{U(q),L(q)}) with $q=0$, we have
\be
  \Ps_{(0)}
  =
  \Us_{(0)}
  =
  (\bfLs^{11}_{(0)})^{-1} \bfLs^{12}_{(0)}
  =
  (\bfLf^{11}_{(0)})^{-1} \bfLf^{12}_{(0)} ,
\Label{U(0)-def}
\ee
where 
$\bfLf_{(0)} 
= \Bf_{(0)} (Dg) \Af^{(0)}$,
according to the definition in 
Eq.~(\ref{Lambda(q)}).
Now, $\bfLf_{(0)}$ admits 
an asymptotic expansion,
$\bfLf_{(0)} 
= \sum_{j=0}^\infty \eps^j \bfLf_{(0,j)}$,
and each of the coefficient matrices $\bfLf_{(0,j)}$
consists of four blocks,
\begin{eqnarray}
    \bfLf^{11}_{(0,j)} &=& \bigg[  \Bf_{(0)}^{12} (D_z g_2)_j
                          + \Bf_{(0)}^{11} (D_z g_1)_{j-1}
                    \bigg] 
                    \Af^{(0)}_{21}, 
\Label{L11(0,q)}\\
    \bfLf^{12}_{(0,j)} &=& \Bf_{(0)}^{12} \left[  (D_y g_2)_j \Af^{(0)}_{12}
                                 + (D_z g_2)_j \Af^{(0)}_{22}
                           \right] \nonumber\\
                           &&\mbox{}+
                  \Bf_{(0)}^{11} \left[  (D_y g_1)_{j-1} \Af^{(0)}_{12}
                                 + (D_z g_1)_{j-1} \Af^{(0)}_{22}
                           \right], 
\Label{L12(0,q)}\\
    \bfLf^{21}_{(0,j)} &=& \Bf_{(0)}^{21} (D_z g_1)_{j-1} \Af^{(0)}_{21}, 
\Label{L21(0,q)}\\
    \bfLf^{22}_{(0,j)} &=&   \Bf_{(0)}^{21} \left[  (D_y g_1)_{j-1} \Af^{(0)}_{12}
                                   + (D_z g_1)_{j-1} \Af^{(0)}_{22}
                             \right].
\Label{L22(0,q)}
\end{eqnarray}
The notation $(\,\cdot\,)_j$ indicates
the $j$th term in the asymptotic expansion
of the quantity inside the parentheses, and
it is understood that such a term is absent
if the subscript is negative.

A direct evaluation shows that
the blocks $\bfLf^{11}_{(0,0)}$ 
and $\bfLf^{12}_{(0,0)}$
are nonzero.
Therefore, $\bfLf^{11}_{(0)}$ and $\bfLf^{12}_{(0)}$
are both $\mathcal{O}(1)$,
and
\be
  \Ps_{(0,0)}
  =
  \Us_{(0,0)}
= (\bfLf^{11}_{(0,0)})^{-1} 
  \bfLf^{12}_{(0,0)}
  =
  \Bf_{(0)}^{12}
  \left[
  \left(D_z g_2\right)_0^{-1}
  \left(D_y g_2\right)_0 \Af^{(0)}_{12} + \Af^{(0)}_{22}
  \right].
\Label{U(0)-aux}
\ee
Here, all the quantities 
are evaluated on $\Mo$,
where the identity
\be
    \left(D_z g_2\right)_0^{-1} \left(D_y g_2\right)_0 = -Dh_0
\ee
holds.
Hence, Eq.~(\ref{U(0)-aux}) implies
\be
  \Ps_{(0,0)}
  =
  \Us_{(0,0)}
  =
  \Bf_{(0)}^{12}
  \left(\Af^{(0)}_{22} - (Dh_0) \Af^{(0)}_{12}\right) .
\Label{U(0,0)}
\ee
Finally, substituting
this expression for
${\tilde P}_{(0,0)}$
into Eq.~(\ref{CSP-cond(1,1)-aux})
and using the identity
$\Af_{21}^{(0)} \Bf^{11}_{(0)} = - \Af_{22}^{(0)} \Bf^{21}_{(0)}$,
we obtain
\be
      (D_z g_2)_0 \psis_{(1,1)} 
    + (D_\eps g_2)_0 
    - (Dh_0) g_{1,0} = 0.
\Label{CSP-cond(1,1)}
\ee
This equation for $\psis_{(1,1)}$ 
is the same as Eq.~(\ref{h(1)}) for $h_1$;
hence, $\psis_{(1,1)}= h_1$ and
$\psis_{(1)} = h_0 + \eps h_1 + \mathcal{O} (\eps^2)$.
This proves the theorem for $q=1$.

\subsection{Proof of Theorem~\ref{1step-Main-thm} for $\mathbf{q=2,3,\ldots}$
\label{ss-1step-ind}}
We prove that
Eqs.~(\ref{U(q-1)-order})--(\ref{psi(q)-exp})
hold for $q+1$,
assuming that they
hold for $0,1,\ldots\,,q$.
By our discussion of
Eqs.~(\ref{Dho-eq})
and (\ref{Dhl-eq}),
Eq.~(\ref{psi(q)-exp})
follows immediately from
Eq.~(\ref{SumU-assum}),
so we need only to consider
Eqs.~(\ref{U(q-1)-order})
and~(\ref{SumU-assum}).

\subsubsection{Establishing Eq.~(\ref{U(q-1)-order})
\label{ss-eq1}}
We first consider Eq.~(\ref{U(q-1)-order}).
The induction hypothesis gives the estimate
$\Us_{(i)}( \cdot\,,\psis_{(i)},\eps)
=\mathcal{O}(\eps^i)$
for $i=0,1,\ldots\,,q-1$.
Also,
$\psis_{(q)} = \psis_{(i)}
+ \mathcal{O}(\eps^{i+1})$
for $i=0,1,\ldots\,,q-1$.
Hence,
$\Us_{(i)}(\cdot\,,\psis_{(q)},\eps)
= \Us_{(i)}(\cdot\,,\psis_{(i)},\eps)
+ \mathcal{O}(\eps^{i+1})$,
from which it follows that
\be
\Us_{(i)}(\cdot\,,\psis_{(q)},\eps)
=  
\mathcal{O}(\eps^i) ,
\qquad
i=0,1,\ldots\,,q-1.
\Label{L-order-aux1}
\ee
In particular,
$\Us_{(0)}(\cdot\,,\psis_{(q)},\eps)
=  \mathcal{O}(1)$,
so
\begin{eqnarray*}
  \Ps_{(q-1)}
=
  \sum_{\ell=0}^{q-1} \Us_{(\ell)}
=
  \mathcal{O}(1)~\mbox{on}~\Ks^{(q)}.
\end{eqnarray*}
This asymptotic estimate
can be used to derive 
asymptotic expansions of 
the blocks of
$\bfLf_{(q)}$.
We begin with
$\bfLs^{11}_{(q)}$.
From Eq.~(\ref{L-expl}),
we have
\be
    \bfLs^{11}_{(q)}
    =
    \bfLf^{11}_{(0)} + \Ps_{(q-1)} \bfLf^{21}_{(0)} .
\Label{L11-expl}
\ee
Since
$\bfLf^{21}_{(0)}
=\mathcal{O}(\eps)$
by
Eq.~(\ref{L21(0,q)}),
we see immediately that
\be
  \bfLs^{11}_{(q)}
= 
    \bfLf^{11}_{(0,0)} 
  + \mathcal{O}(\eps) .
\Label{L11-order}
\ee
Next, we examine
the block $\bfLs^{12}_{(q)}$.
From Eq.~(\ref{L12-rec}), we have
\[
  \bfLs^{12}_{(q)} 
=
    \Us_{(q-1)} \bfLs^{22}_{(q-1)}
  - \Us_{(q-1)} \bfLs^{21}_{(q-1)} \Us_{(q-1)} 
  + \frac{d\Us_{(q-1)}}{dt} .
\]
First,
$\Us_{(q-1)}(\cdot\,,\psis_{(q)},\eps)
= \mathcal{O}(\eps^{q-1})$
by Eq.~(\ref{L-order-aux1}).
Also,
$\bfLs^{21}_{(q-1)}
= \bfLf^{21}_{(0)}
= \mathcal{O}(\eps)$
on $\Ks^{(q)}$
by Eqs.~(\ref{L-expl})
and~(\ref{L21(0,q)}).
Moreover,
$\bfLs^{22}_{(q-1)}
= \bfLf^{22}_{(0)}
- \bfLf^{21}_{(0)}
  \Ps_{(q-2)}
= \mathcal{O}(\eps)$
by Eqs.~(\ref{L-expl})
and~(\ref{L22(0,q)}).
Finally,
by applying Lemma~\ref{DUg-lem}
with $V=\Us_{(q-1)}$,
we find that $d\Us_{(q-1)}/dt$
is $\mathcal{O}(\eps^q)$.
Putting these estimates together,
we obtain the estimate
\be
  \bfLs^{12}_{(q)}
=
    \eps^q \bfLs^{12}_{(q,q)}
  + \mathcal{O}(\eps^{q+1}) ,
\Label{L12-order}
\ee
where we grouped all of
the $\mathcal{O}(\eps^q)$ terms
into $\eps^q \bfLs^{12}_{(q,q)}$.
By combining
the definition~(\ref{U(q),L(q)})
with Eqs.~(\ref{L11-order})
and~(\ref{L12-order}),
we derive the desired estimate,
$\Us_{(q)}
= (\bfLs^{11}_{(q)})^{-1}
  \bfLs^{12}_{(q)}
= \mathcal{O}(\eps^q)$.

\paragraph{Remark.}
While the estimates of
$\bfLs^{21}_{(q)}$
and
$\bfLs^{22}_{(q)}$
are not needed here,
they will be needed in
Section~\ref{s-Full}.
First, $\bfLs^{21}_{(q)}
= \bfLf^{21}_{(0)}
= \mathcal{O}(\eps)$
on $\Ks^{(q)}$,
by Eqs.~(\ref{L-expl})
and~(\ref{L21(0,q)}).
Then,
$\bfLs^{22}_{(q)}
= \bfLf^{22}_{(0)}
- \bfLf^{21}_{(0)}
  \Ps_{(q-1)}$
by Eq.~(\ref{L-expl}).
Now, $\bfLf^{22}_{(0)}
= \mathcal{O}(\eps)$
by Eq.~(\ref{L22(0,q)}),
and thus
the discussion for
the size of $\bfLs^{11}_{(q)}$
also yields that
$\bfLs^{22}_{(q)}
=\mathcal{O}(\eps)$.
Putting the estimates of
this section together,
we obtain
\begin{eqnarray}
    \bfLs_{(q)}( \cdot\,,\psis_{(q)},\eps)
    =
    \left(
          \begin{array}{cc}
          \bfLf^{11}_{(0,0)} + \mathcal{O}(\eps)
&         \eps^q\bfLs^{12}_{(q,q)} + \mathcal{O}(\eps^{q+1})\\
          \eps\bfLf^{21}_{(0,1)} + \mathcal{O}(\eps^2)
&         \eps\bfLs^{22}_{(q,1)} + \mathcal{O}(\eps^2)
          \end{array}
    \right) .
\Label{L-order}
\end{eqnarray}

\subsubsection{Establishing Eq.~(\ref{SumU-assum})
\label{ss-eq2}}
Next, we consider Eq.~(\ref{SumU-assum}).
The induction hypothesis gives the estimate
$\Af_{21}^{(0)}
[\Ps_{(i)}(\cdot\,,\psis_{(i)},\eps)
 \Bf^{21}_{(0)} + \Bf^{11}_{(0)}]
= - \sum_{j=0}^i \eps^j Dh_j
+ \mathcal{O}(\eps^{i+1})$
for $i=0,1,\ldots\,,q-1$.
Our goal is to show that this equation 
also holds for $i=q$.
We first show that
the terms up to and including
$\mathcal{O}(\eps^{q-1})$
in both members of the equation 
agree for $i=q$.
Then we analyze the
terms of $\mathcal{O}(\eps^q)$.

By the induction hypothesis,
we have the asymptotic expansion
\begin{eqnarray}
    \Af^{(0)}_{21} 
    \left[\Ps_{(q-1)}(\cdot\,,\psis_{(q-1)},\eps) \Bf_{(0)}^{21} 
          + \Bf_{(0)}^{11}\right] 
    &=& 
    - \sum_{j=0}^{q-1} \eps^j Dh_j 
    + \mathcal{O}(\eps^q) .
\end{eqnarray}
Also by the induction hypothesis,
$\psis_{(q)} = \psis_{(q-1)}
+\mathcal{O}(\eps^q)$.
Hence,
\begin{eqnarray}
    \Af^{(0)}_{21} 
    \left[\Ps_{(q-1)}(\cdot\,,\psis_{(q)},\eps) \Bf_{(0)}^{21} 
          + \Bf_{(0)}^{11}\right] 
    &=& 
    - \sum_{j=0}^{q-1} \eps^j Dh_j 
    + \mathcal{O}(\eps^q) .
\Label{P-aux}
\end{eqnarray}
The definition~(\ref{expl-aux1}) of $\Ps_{(q)}$
yields the update formula
\be
    \Ps_{(q)}
=
    \Ps_{(q-1)} + \Us_{(q)} .
\Label{P-rec}
\ee
We already showed that
$\Us_{(q)}(\cdot\,,\psis_{(q)},\eps)
= \mathcal{O}(\eps^q)$,
so Eq.~(\ref{P-rec})
implies that the asymptotic expansions of
$\Ps_{(q)}(\cdot\,,\psis_{(q)},\eps)$
and
$\Ps_{(q-1)}(\cdot\,,\psis_{(q)},\eps)$
agree up to and including
terms of $\mathcal{O}(\eps^{q-1})$.
The same, then, holds for
the asymptotic expansions of
$\Af^{(0)}_{21}
[\Ps_{(q)}(\cdot\,,\psis_{(q)},\eps) \Bf_{(0)}^{21} 
      + \Bf_{(0)}^{11}]$
and
$\Af^{(0)}_{21} 
[\Ps_{(q-1)}(\cdot\,,\psis_{(q)},\eps) \Bf_{(0)}^{21} 
      + \Bf_{(0)}^{11}]$.
Therefore,
\begin{eqnarray}
    \Af^{(0)}_{21} 
    \left[\Ps_{(q)}(\cdot\,,\psis_{(q)},\eps) \Bf_{(0)}^{21} 
          + \Bf_{(0)}^{11}\right] 
    &=& 
    - \sum_{j=0}^{q-1} \eps^j Dh_j 
    + \mathcal{O}(\eps^q) .
\end{eqnarray}
In other words,
\begin{eqnarray}
    \Af^{(0)}_{21} 
    \left[\Ps_{(q,0)} \Bf_{(0)}^{21} + \Bf_{(0)}^{11}\right] 
    &=& 
    - Dh_0 ,
\Label{AP(q,0)B=-Dho}\\
    \Af^{(0)}_{21} 
    \Ps_{(q,j)} \Bf_{(0)}^{21} 
    &=& 
    - Dh_j,
    \hskip0.1truein \mbox{for all} 
    \hskip0.1truein j=1,\ldots\,,q-1 ,
\Label{AP(q,j)B=-Dhj}
\end{eqnarray}
which establishes 
Eq.~(\ref{SumU-assum})
for all terms up to and including
$\mathcal{O}(\eps^{q-1})$.

It remains to show that
the terms of $\mathcal{O}(\eps^q)$
in both members of Eq.~(\ref{SumU-assum})
agree,
that is,
\begin{eqnarray}
    \Af^{(0)}_{21} \Ps_{(q,q)} \Bf_{(0)}^{21} 
    &=& 
    - Dh_q .
\Label{AP(q,q)B=-Dhq}
\end{eqnarray}
We achieve this
by deriving an explicit formula for
$\Af^{(0)}_{21} \Ps_{(q,q)} \Bf_{(0)}^{21}$
and comparing it to that for $Dh_q$,
which is given in the Appendix
(Lemma~\ref{Dhq-lem}).
We proceed in two steps.
In step one, 
we express $\Af^{(0)}_{21} \Ps_{(q,q)} \Bf_{(0)}^{21}$ 
in terms of
$\Ps_{(q-1,0)},\ldots\,,\Ps_{(q-1,q-1)}$.
Then, in step two, 
we obtain the explicit formula for 
$\Af^{(0)}_{21} \Ps_{(q,q)} \Bf_{(0)}^{21}$
in terms of the vector field 
and of $Dh_i$, $i=0,1,\ldots\,,q-1$.

Step 1.~~Recall the update formula~(\ref{P-rec}),
$\Ps_{(q)} = \Ps_{(q-1)} + \Us_{(q)}$.
Using the definition~(\ref{U(q),L(q)}) 
of $\Us_{(q)}$
and the explicit formula~(\ref{L-expl})
for $\bfLs_{(q)}$,
we can express $\Us_{(q)}$ in terms of
$\bfLf_{(0)}$ and $\Ps_{(q-1)}$.
In particular,
Eq.~(\ref{L-order}) implies that
$\Us_{(q,q)}
=(\bfLf^{11}_{(0,0)})^{-1}
\bfLs^{12}_{(q,q)}$.
Also,
Eq.~(\ref{L-expl})
gives
\be
    \bfLs^{12}_{(q)}
    =
    \bfLf^{12}_{(0)} - \bfLf^{11}_{(0)} \Ps_{(q-1)}
    + \Ps_{(q-1)} \bfLf^{22}_{(0)}
    - \Ps_{(q-1)} \bfLf^{21}_{(0)} \Ps_{(q-1)}
    + \frac{d\Ps_{(q-1)}}{dt} .
\Label{L12-expl}~~
\ee
It follows that
\begin{eqnarray}
    \Us_{(q,q)}
    &=&
    (\bfLf^{11}_{(0,0)})^{-1}
    \left[
     \bfLf^{12}_{(0,q)}
     - \left(\bfLf^{11}_{(0)} \Ps_{(q-1)}\right)_q
     + \left(\Ps_{(q-1)} \bfLf^{22}_{(0)}\right)_q  \right.\nonumber\\
     &&\mbox{}\left.
     - \left(\Ps_{(q-1)} \bfLf^{21}_{(0)} \Ps_{(q-1)}\right)_q
     + \left(\frac{d\Ps_{(q-1)}}{dt}\right)_q 
    \right] ,
\Label{U(q,q)-aux1}
\end{eqnarray}
where
we recall
the notational convention
that
$(\,\cdot\,)_q$
stands for 
the coefficient of
the
$\mathcal{O}(\eps^q)$
term in
the asymptotic expansion of
the quantity in parentheses.
Using Lemma~\ref{DUg-lem}
with $V=\Ps_{(q-1)}$
and the fact that
$\bfLf^{22}_{(0,0)}$
and
$\bfLf^{21}_{(0,0)}$
are both zero,
we rewrite
Eq.~(\ref{U(q,q)-aux1})
as
\begin{eqnarray}
    \Us_{(q,q)}
    &=&
    (\bfLf^{11}_{(0,0)})^{-1}
    \left[
J_1 + ( J_2 - \Lambda^{11}_{(0,0)} \Ps_{(q-1,q)} ) + J_3 + J_4 + J_5
   \right],~~~~~
\Label{U(q,q)}
\end{eqnarray}
where
\begin{eqnarray}
    J_1
=
    \bfLf^{12}_{(0,q)},~~~~~
    J_2
=
     - \sum_{\ell=0}^{q-1} 
     \bfLf^{11}_{(0,q-\ell)} \Ps_{(q-1,\ell)} ,~~~~~
    J_3
=
    \sum_{\ell=0}^{q-1} 
     \Ps_{(q-1,\ell)} \bfLf^{22}_{(0,q-\ell)} , ~~~~~\nonumber\\
    J_4
=
    - \sum_{i=0}^{q-1}\sum_{j=0}^{q-1-i}
            \Ps_{(q-1,j)} \bfLf^{21}_{(0,q-i-j)} \Ps_{(q-1,i)} ,~~~~~
    J_5
=
    \sum_{\ell=0}^{q-1} 
            \frac{d\Ps_{(q-1,\ell)}}{dy} g_{1,q-1-\ell} .~~~~~
\end{eqnarray}
Substituting the expression~(\ref{U(q,q)})
into the update formula~(\ref{P-rec}) for $\Ps_{(q)}$,
we find
\begin{eqnarray}
    \Af^{(0)}_{21} 
    \Ps_{(q,q)} 
    \Bf_{(0)}^{21}
    &=&
    \Af^{(0)}_{21} 
    (\bfLf^{11}_{(0,0)})^{-1}
    \left[
          J_1 + J_2 + J_3
          + J_4 + J_5
   \right] 
    \Bf_{(0)}^{21} ,
\Label{P(q,q)}
\end{eqnarray}
Step 2.~~We rewrite
the terms $J_1,\ldots\,,J_5$
by means of the induction hypothesis and
the explicit formulas~(\ref{L11(0,q)})--(\ref{L22(0,q)})
for the blocks of $\bfLf_{(0)}$.

Equation~(\ref{L11(0,q)})
and the identity
$\Af^{(0)}_{21}\Bf_{(0)}^{12}
= I_n$
imply that
\be
    \Af^{(0)}_{21}
    \left(\bfLf^{11}_{(0,0)}\right)^{-1} 
    = 
    \left((D_z g_2)_0\right)^{-1} 
    \Af^{(0)}_{21}.
\Label{L11-0}
\ee
Here,
$(D_z g_2)_0$
stands for
the leading order term in
the asymptotic expansion of
$(D_z g_2)(\cdot\,,\psis_{(q)},\eps)$.
Since
$\psis_{(q)}$
and
$h_\eps$
agree 
up to and including
$\mathcal{O}(\eps^q)$
terms by assumption,
the asymptotic expansions of
$(D_z g_2)(\cdot\,,\psis_{(q)},\eps)$
and
$(D_z g_2)(\cdot\,,h_\eps,\eps)$
also agree
up to and including
$\mathcal{O}(\eps^q)$
terms.
For the remainder of
this section,
it does not matter
whether
quantities are evaluated on
$\Ks^{(q)}$
or on
$\Me$,
since only the coefficients of
$\eps^q$ or lower 
appear in our formulas.
Accordingly, we make
no distinction between
the asymptotic expansions
of a quantity 
evaluated on
the two manifolds.

Using
Eq.~(\ref{L12(0,q)})
and
the identities
$\Bf_{(0)}^{12}=(\Af^{(0)}_{21})^{-1}$,
$\Bf_{(0)}^{21}=(\Af^{(0)}_{12})^{-1}$,
and
$\Bf_{(0)}^{11}=-\Bf_{(0)}^{12}\Af^{(0)}_{22}\Bf_{(0)}^{21}$,
we find
\begin{eqnarray}
\lefteqn{
    \Af^{(0)}_{21} J_1 \Bf_{(0)}^{21}
    =
      (D_y g_2)_q
    + (D_z g_2)_q \Af^{(0)}_{22} \Bf_{(0)}^{21}}&& \nonumber\\
    &&\mbox{}-
      \Af^{(0)}_{22} \Bf_{(0)}^{21} (D_y g_1)_{q-1}
    \mbox{}
    - \Af^{(0)}_{22} \Bf_{(0)}^{21} (D_z g_1)_{q-1} \Af^{(0)}_{22} \Bf_{(0)}^{21}.
\Label{AL12B}
\end{eqnarray} 
Next, substituting for
$\Af^{(0)}_{21} \Ps_{(q-1,\ell)}\Bf_{(0)}^{21}$
from
the induction hypothesis~(\ref{SumU-assum}),
we obtain
\begin{eqnarray}
\lefteqn{
     \Af^{(0)}_{21}
     J_2
     \Bf_{(0)}^{21}
    =
    \sum_{\ell=0}^{q-1} (D_z g_2)_{q-\ell} Dh_\ell 
    - (D_z g_2)_q 
    \Af^{(0)}_{22} \Bf_{(0)}^{21}}&&
    \nonumber\\ 
    &&\mbox{}
    - \sum_{\ell=0}^{q-1} 
    \Af^{(0)}_{22} \Bf_{(0)}^{21} 
    (D_z g_1)_{q-1-\ell} Dh_\ell
    + \Af^{(0)}_{22} \Bf_{(0)}^{21} 
    (D_z g_1)_{q-1} 
    \Af^{(0)}_{22} \Bf_{(0)}^{21}.
\Label{AL11PB}
\end{eqnarray}
Then, using
Eq.~(\ref{L22(0,q)})
and the assumptions of
the lemma,
we find
\begin{eqnarray}
\lefteqn{
  \Af^{(0)}_{21}
  J_3
  \Bf_{(0)}^{21}
    =\mbox{}
    - \sum_{\ell=0}^{q-1} 
    Dh_\ell (D_y g_1)_{q-1-\ell}
    - \sum_{\ell=0}^{q-1} 
    Dh_\ell (D_z g_1)_{q-1-\ell} 
    \Af^{(0)}_{22} \Bf_{(0)}^{21}} \hspace{4em}&&
    \nonumber\\
    &&\mbox{}
    + \Af^{(0)}_{22} \Bf_{(0)}^{21} 
    (D_y g_1)_{q-1} 
    + \Af^{(0)}_{22} \Bf_{(0)}^{21} 
    (D_z g_1)_{q-1} 
    \Af^{(0)}_{22} \Bf_{(0)}^{21} .
\Label{APL22B}
\end{eqnarray}
In the same vein,
we use 
the induction hypothesis
on $J_4$,
\begin{eqnarray}
\lefteqn{\hspace{-4em}
  \Af^{(0)}_{21}
  J_4 
  \Bf_{(0)}^{21}
    =\mbox{}
    - \sum_{i=0}^{q-1}\sum_{j=0}^{q-1-i} 
    Dh_j (D_z g_1)_{q-1-i-j} Dh_i
    + \sum_{i=0}^{q-1}
    \Af^{(0)}_{22} \Bf_{(0)}^{21} 
    (D_z g_1)_{q-1-i} Dh_i }
    \nonumber\\
    &&\hspace{-2em}\mbox{}
    + \sum_{j=0}^{q-1}\
    Dh_j (D_z g_1)_{q-1-j} 
    \Af^{(0)}_{22} \Bf_{(0)}^{21}
    - \Af^{(0)}_{22} \Bf_{(0)}^{21} 
    (D_z g_1)_{q-1} 
    \Af^{(0)}_{22} \Bf_{(0)}^{21} .  
\Label{APL21PB}
\end{eqnarray}
The terms in Eq.~(\ref{P(q,q)})
containing $A_{22}^{(0)}$
sum to zero,
which may be seen as follows.
The second and fourth
terms in
(\ref{AL12B})
cancel against
the second and fourth
terms in
(\ref{AL11PB});
the third term in
(\ref{AL12B})
cancels against 
the third term in
(\ref{APL22B});
the third term in
(\ref{AL11PB})
cancels against 
the second term in
(\ref{APL21PB}); and
the second and fourth
terms in
(\ref{APL22B})
cancel against the
third and fourth
terms in
(\ref{APL21PB}).
These cancellations
were to be expected
because 
the approximation
should be independent of
the choice of
$\Af^{(0)}$.

Carrying out
the same type of calculation 
as above, we obtain
\begin{eqnarray}
    \Af^{(0)}_{21}
    J_5
    \Bf_{(0)}^{21}
    &=&
    - \sum_{\ell=0}^{q-1} 
    D^2h_\ell g_{1,q-1-\ell} ,
\Label{DPg}
\end{eqnarray}
where
we have used 
the symmetry of
the bilinear form
$D^2h_\ell$.

Equations~(\ref{L11-0})--(\ref{DPg}),
together with
the observed cancellations,
yield
\begin{eqnarray}
\lefteqn{
    \Af^{(0)}_{21} \Ps_{(q,q)} \Bf_{(0)}^{21} } \nonumber \\
    &=&
    \left((D_z g_2)_0\right)^{-1}
    \left[
    (D_y g_2)_q %\,\fbox{\fbox{U4}}
    + \sum_{\ell=0}^{q-1}
      (D_z g_2)_{q-\ell} Dh_\ell 
    \right.%\nonumber\\
    - \sum_{\ell=0}^{q-1} D^2h_\ell g_{1,q-1-\ell} \nonumber\\
    &&\mbox{}
    - \sum_{\ell=0}^{q-1} Dh_\ell (D_y g_1)_{q-1-\ell} 
    \left.
    - \sum_{i=0}^{q-1}\sum_{j=0}^{q-1-i} 
    Dh_j (D_z g_1)_{q-1-i-j} Dh_i
    \right] .
\Label{AP(q,q)B-exp}
\end{eqnarray}
A term-by-term comparison
with the expression for $-Dh_q$
given in the Appendix, Eq.~(\ref{Dhq}),
shows that
$\Af_{21}^{(0)} \Ps_{(q,q)} \Af^{21}_{(0)}
= - Dh_q$.
Thus, the proof of
Theorem~\ref{1step-Main-thm}
is complete.

\paragraph{Remark.}
In general, the error term is nontrivial,
as can already be seen at $q=0$.
The equation determining $\psis_{(0,1)}$ is
\be
      (D_z g_2)_0 \psis_{(0,1)} 
    + (D_\eps g_2)_0
    - \Af^{(0)}_{22} \Bf_{(0)}^{21} g_{1,0} = 0.
\Label{CSP-cond(0,1)}
\ee
This equation is
not the same as Eq.~(\ref{h(1)}),
which determines $h_1$.
Where Eq.~(\ref{h(1)})
has the term $Dh_0$,
Eq.~(\ref{CSP-cond(0,1)})
has the term
$\Af^{(0)}_{22}\Bf_{(0)}^{21}$.
When the slow manifold is nonlinear,
$Dh_0$
depends on 
$y$,
whereas
$\Af^{(0)}_{22}\Bf_{(0)}^{21}$
is a constant matrix.
Therefore, in general
$\psis_{(0,1)} \neq h_1$,
and the strongest claim
we can make is
$\psis_{(0)} = h_0 + \mathcal{O} (\eps)$.
A similar argument applies to higher values of $q$.

%%%%%%%%%%%%%%%%%%%%%%%%%%%%%%%%%%%%%%%%%%%%%%%%%%%%%%%%%%%%%%%%%%%%%%%
%%%%%%%%%%%%%%%%%%%%%%%%%%%%%%%%%%%%%%%%%%%%%%%%%%%%%%%%%%%%%%%%%%%%%%%
%%%%%%%%%%%%%%%%%%%%%%%%%%%%%%%%%%%%%%%%%%%%%%%%%%%%%%%%%%%%%%%%%%%%%%%

\section{Analysis of the Full CSP Method\label{s-Full}}
\setcounter{equation}{0}
We now return to
the full CSP method
and prove Theorem~\ref{t-CSPMq}.
Since the full CSP method
and the one-step CSP method
start from the same basis, 
the conditions~(\ref{CSPM-q})
and~(\ref{CSPM-q-1step})
are the same for $q=0$,
\be
    \Bf^1_{(0)} g = 0 .
\ee
Therefore,
we can choose
$\psif_{(0)} = \psis_{(0)} = h_0$.

%%%%%%%%%%%%%%%%%%%%%%%%%%%%%%%%%%%%%%%%%%%%%%%%%%%%%%%%%%%%%%%%%%%%%%%
%%%%%%%%%%%%%%%%%%%%%%%%%%%%%%%%%%%%%%%%%%%%%%%%%%%%%%%%%%%%%%%%%%%%%%%
%%%%%%%%%%%%%%%%%%%%%%%%%%%%%%%%%%%%%%%%%%%%%%%%%%%%%%%%%%%%%%%%%%%%%%%

\subsection{Proof of Theorem~\ref{t-CSPMq} for $\mathbf{q=1}$
\label{s-full-q=1}}

In this section,
we carry out
the first iteration of
the full CSP method
and
determine the 
resulting approximation
$\Kf^{(1)}$
of the slow manifold.
We then compare
$\Kf^{(1)}$
and $\Ks^{(1)}$.

The update quantities
$\Uf_{(0)}$ 
and 
$\Lf_{(0)}$
follow from the definition~(\ref{U(q),L(q)}),
\begin{eqnarray}
    \Uf_{(0)} 
=
    (\bfLf^{11}_{(0)})^{-1} \bfLf^{12}_{(0)} ,
\qquad
    \Lf_{(0)} 
=
    \bfLf^{21}_{(0)} (\bfLf^{11}_{(0)})^{-1} .
\Label{Uf(0),Lf(0)-def}
\end{eqnarray}
(We recall that we use the same notation
$\Uf_{(0)}$ and $\Lf_{(0)}$
for the full matrix and
the nonzero block.)
In particular,
Eqs.~(\ref{Uf(0),Lf(0)-def})
and
(\ref{U(0)-def})
imply that
$\Uf_{(0)}= \Us_{(0)}$.
Next, we update the matrix $\Bf_{(0)}$.
Following Eq.~(\ref{B(q+1)}),
we find
\be
    \Bf_{(1)} = (I-\Lf_{(0)})(I+\Uf_{(0)}) \Bf_{(0)}.
\Label{i=1-aux}
\ee
The upper and lower
row blocks of $\Bf_{(1)}$ are
\begin{eqnarray}
    \Bf^1_{(1)}
    &=&
      \Bf^1_{(0)}
    + \Uf_{(0)} \Bf^2_{(0)},
\Label{Bf1(1)-expl}\\
    \Bf^2_{(1)}
    &=&
       (I - \Lf_{(0)}\Uf_{(0)})
       \Bf^2_{(0)}
    - \Lf_{(0)} \Bf^1_{(0)}.
\Label{Bf2(1)-expl}
\end{eqnarray}
Since
$\Ps_{(0)}=\Us_{(0)}=\Uf_{(0)}$
and
$\psif_{(0)} = \psis_{(0)}$,
Eqs.~(\ref{B-expl})
and~(\ref{Bf1(1)-expl})
imply that
\be
    \Bf^1_{(1)} = \Bs^1_{(1)} ,
\ee
so
after the first iteration
the CSP condition
is the same as for
the one-step method.
Therefore,
$\psif_{(1)}=\psis_{(1)}$
and, by Theorem~\ref{1step-Main-thm},
\be
    \psif_{(1)} = h_0 + \eps h_1 + {\cal O}(\eps^2) .
\Label{psif(1)-exp}
\ee
This proves Theorem~\ref{t-CSPMq} for $q=1$.

%%%%%%%%%%%%%%%%%%%%%%%%%%%%%%%%%%%%%%%%%%%%%%%%%%%%%%%%%%%%%%%%%%%%%%%
%%%%%%%%%%%%%%%%%%%%%%%%%%%%%%%%%%%%%%%%%%%%%%%%%%%%%%%%%%%%%%%%%%%%%%%
%%%%%%%%%%%%%%%%%%%%%%%%%%%%%%%%%%%%%%%%%%%%%%%%%%%%%%%%%%%%%%%%%%%%%%%

\subsection{The Induction Hypothesis
\label{ss-der-full}}
So far,
we have established the identities
$\Bf^1_{(0)}=\Bs^1_{(0)}$
and
$\Bf^1_{(1)}=\Bs^1_{(1)}$,
from which we could conclude that
$\Kf^{(0)} = \Ks^{(0)}$
and
$\Kf^{(1)} = \Ks^{(1)}$.
In general, though,
it is not true that
$\Bf^1_{(q)} = \Bs^1_{(q)}$
for higher values of $q$,
as we now demonstrate.

In the one-step CSP method,
Eq.~(\ref{B-expl}) yields
\[
    \Bs^1_{(2)}
=
      B^1_{(0)}
    + (\Us_{(0)}+\Us_{(1)})
      B^2_{(0)} .
\]
By contrast, in the full CSP method,
we obtain from Eq.~(\ref{B(q+1)})
\be
    \Bf^1_{(2)}
=
      \left(  I_n
            - \Uf_{(1)}\Lf_{(0)}
      \right)
      B^1_{(0)}
    + \left(  \Uf_{(0)}
            + \Uf_{(1)}
            - \Uf_{(1)}\Lf_{(0)}\Uf_{(0)}
      \right)
      B^2_{(0)} .
\Label{B(2)}
\ee
The rows of
$B^1_{(0)}$
and
$B^2_{(0)}$
are linearly independent,
as can be seen from
Eq.~(\ref{B(0)}),
so the presence of
the premultiplier of
$B^1_{(0)}$
in the expression~(\ref{B(2)})
implies that
$\Bf^1_{(2)}\neq \Bs^1_{(2)}$.
A similar argument shows that
$\Bf^1_{(q)}\neq \Bs^1_{(q)}$
for $q = 2, 3, \ldots\,$.
Consequently, the proof of 
Theorem~\ref{t-CSPMq} for $q=1$
given in Section~\ref{s-full-q=1}
does not generalize
to higher values of $q$.

The matrix $\Bs^1_{(q)}$ has
an important property.
Using
Eq.~(\ref{B-expl}),
we write
\[
    \Bs^1_{(q)} 
    = 
    \left(\Ps_{(q-1)} B_{(0)}^{21} + B_{(0)}^{11}, \,
          B_{(0)}^{12} \right) 
    = 
    B_{(0)}^{12}
    \left(A^{(0)}_{21}
          \left[  \Ps_{(q-1)} B_{(0)}^{21}
          + B_{(0)}^{11}
          \right], \,
            I_n
    \right) .
\]
Given the induction hypothesis~(\ref{SumU-assum}),
we rewrite this expression once more,
\be
    \Bs^1_{(q)}
    = 
    B_{(0)}^{12}
    \left(- \sum_{j=0}^{q-1} 
            \eps^j Dh_j
          + \mathcal{O}(\eps^q) , \,
            I_n
    \right) .
\Label{B1-prop}
\ee
Take any $y\in K$,
and let the points
$\tilde{Q} \in \Ks^{(q-1)}$,
$Q \in \Kf^{(q-1)}$, and
$Q' \in \Me$
be defined by
\[
\tilde{Q} = (y,\psis_{(q-1)}(y,\eps)) , \quad
Q = (y,\psif_{(q-1)}(y,\eps)) , \quad
Q' = (y,h_\eps(y)) .
\]
The $n$ row vectors of the matrix
$( - Dh_\eps(y), \, I_n)$
form an exact basis for
${\cal N}_{Q'}\Me$,
the space normal to $\Me$ at $Q'$.
Therefore,
by Eq.~(\ref{B1-prop}),
$\Bs^1_{(q)}(\tilde{Q})$
is
a linear combination of
the basis vectors of
${\cal N}_{Q'}\Me$,
up to and including 
terms of
$\mathcal{O}(\eps^{q-1})$,
via the invertible matrix 
$B_{(0)}^{12}$.
Hence,
the columns of
$\Bs^1_{(q)}(\tilde{Q})$
form a basis for
${\cal N}_{Q'}\Me$
up to and including
terms of $\mathcal{O}(\eps^{q-1})$.
This property of
$\Bs^1_{(q)}(\tilde{Q})$
was central to the proof of
Theorem~\ref{1step-Main-thm}.
We seek to prove
a similar result for
the rows of
$\Bf^1_{(q)}(Q)$.

The rows of
$\Bf_{(q)}(Q)$
can be written as
linear combinations of
the rows of
$\Bs_{(q)}(\tilde{Q})$,
\be
  \Bf_{(q)}(Q) 
= 
  T_{(q)}(y,\eps) 
  \Bs_{(q)} (\tilde{Q}) ,
\Label{Bf=TB}
\ee
because $\Bs_{(q)}(\tilde{Q})$
is invertible
(see Eq.~(\ref{B-expl})).
In terms of
the constituent blocks,
\begin{eqnarray}
    \Bf^1_{(q)}
    &=&
      T^{11}_{(q)} \Bs^1_{(q)}
    + T^{12}_{(q)} \Bs^2_{(q)},
\Label{Bf1=TB}\\
    \Bf^2_{(q)}
    &=&
      T^{21}_{(q)} \Bs^1_{(q)}
    + T^{22}_{(q)} \Bs^2_{(q)}.
\Label{Bf2=TB}
\end{eqnarray}
Equation~(\ref{Bf1=TB}) shows
that the requirement that
the rows of
$\Bf^1_{(q)}(Q)$
span
${\cal N}_{Q'}\Me$
up to and including
terms of
${\cal O}(\eps^{q-1})$
is equivalent to
the conditions
\be
  T^{11}_{(q)}(y,\eps)
=
  \mathcal{O}(1)
  \mbox{ and invertible} , \quad
  T^{12}_{(q)}(y,\eps)
=
  {\cal O}(\eps^q) .
\Label{T12,T11-assum}
\ee
Assume for the moment that
these conditions are satisfied.
Then the CSP condition~(\ref{CSPM-q})
after the $q$th iteration
can be recast as
\[
  \left[ T^{11}_{(q)} (y,\eps)
  \Bs^1_{(q)} (y,\psif_{(q-1)} (y, \eps), \eps)
    +
  T^{12}_{(q)} (y,\eps)
  \Bs^2_{(q)} (y,\psif_{(q-1)}(y, \eps), \eps) \right]
  g(y,z,\eps)
    = 0
\]
or, since $T^{11}_{(q)}(y,\eps)$ is invertible,
\be
  \Bs^1_{(q)}
  g
    +
  \left(T^{11}_{(q)}\right)^{-1} 
  T^{12}_{(q)}
  \Bs^2_{(q)}
  g
    =
    0.
\Label{CSPf-cond(q)}
\ee
The second term
is at least of ${\cal O}(\eps^q)$,
by the second assumption in
Eq.~(\ref{T12,T11-assum}),
so the terms of
${\cal O}(\eps^j)$
in
Eqs.~(\ref{CSPM-q-1step})
and~(\ref{CSPf-cond(q)})
are equal
for
$j=0,1,\ldots\,,q-1$.
At ${\mathcal O}(\eps^q)$,
the two equations 
differ by the term
$(T^{11}_{(q,0)})^{-1}
T^{12}_{(q,q)}
\Bs^2_{(q,0)}
g(y,\psif_{(q,0)},\eps)$.
Since
the ${\mathcal O}(1)$ terms of 
the two equations agree, 
it follows that
$\psif_{(q,0)}=\psis_{(q,0)}=h_0$
and, therefore,
$g(y,\psif_{(q,0)},\eps) = 0$.
Hence,
Eqs.~(\ref{CSP-cond(q)})
and~(\ref{CSPf-cond(q)})
agree
up to and including
terms of
${\cal O}(\eps^{q})$,
so
Eq.~(\ref{CSPf-cond(q)})
produces
the asymptotic expansion of
the slow manifold
up to and including
terms of
${\cal O}(\eps^{q})$,
by Theorem
\ref{1step-Main-thm}.

To complete the proof of Theorem~\ref{t-CSPMq},
we need to verify
the conditions~(\ref{T12,T11-assum})
for $q=2, 3, \ldots\,$,
which we do by
induction on $q$.
The induction hypothesis is
\begin{eqnarray}
  T_{(q)}(\cdot\,,\psif_{(q-1)},\eps)
  &=&
  \left(\begin{array}{cc}
          I_n + {\cal O}(\eps^2)
        & \eps^q T^{12}_{(q,q)}
        + {\cal O}(\eps^{q+1})\\
          \eps T^{21}_{(q,1)}
        + {\cal O}(\eps^2)
        & I_m + {\cal O}(\eps^2)         
        \end{array}
  \right) ,
\Label{T(q)-exp}\\
   \psif_{(q)} (\cdot\,,\eps)
   &=& 
    \sum_{j=0}^q \eps^j h_j + \mathcal{O} (\eps^{q+1}) .
\Label{psif(q)-exp}
\end{eqnarray}

\subsection{Proof of Theorem~\ref{t-CSPMq}
for $\mathbf{q=2,3,\ldots}$
\label{ss-full-ind}}
In this section,
we carry out 
the induction step
of the proof.
We assume that
Eqs.~(\ref{T(q)-exp})
and (\ref{psif(q)-exp})
hold for
$0,1,\ldots\,,q$
and prove that 
they also hold for
$q+1$.
It suffices to establish
Eq.~(\ref{T(q)-exp});
Eq.~(\ref{psif(q)-exp})
follows immediately from
Eq.~(\ref{T(q)-exp})
and our discussion of
the CSP condition~(\ref{CSPf-cond(q)}).

Before carrying out the induction step,
we derive an update formula for
$T_{(q)}$.
Using Eq.~(\ref{Bf=TB})
with $q$ replaced by $q+1$,
we obtain
\be
    T_{(q+1)}
    =
    \Bf_{(q+1)}
    \As^{(q+1)} .
\Label{T(q+1)}
\ee
(Here, we used the identity
$(\Bs_{(q+1)})^{-1}=\As^{(q+1)}$.)
Next, we use
the update formulas~(\ref{B(q+1)})
and~(\ref{A-rec})
for
$\Bf_{(q+1)}$ and
$\As^{(q+1)}$,
respectively,
to rewrite
Eq.~(\ref{T(q+1)}),
\be
    T_{(q+1)}
    =
    \left(I-\Lf_{(q)}\right)
    \left(I+\Uf_{(q)}\right)
    T_{(q)}
    (I-\Us_{(q)}) .
\Label{T-rec}
\ee
Equation~(\ref{Bf=TB})
also relates
$\Af^{(q+1)}$ to $\As^{(q+1)}$,
\be
  \Af^{(q+1)} 
= 
  \As^{(q+1)} 
  \left(T_{(q+1)}\right)^{-1} .
\Label{Af=ATinv}
\ee
Taking
$C=\left(T_{(q)}\right)^{-1}$
in Eq.~(\ref{L-transf}),
we express $\bfLf_{(q)}$
in terms of
$\bfLs_{(q)}$,
\begin{eqnarray*}
    \bfLf_{(q)} =   T_{(q)} \bfLs_{(q)} \left(T_{(q)}\right)^{-1}
                  - T_{(q)} \frac{d\left(T_{(q)}\right)^{-1}}{dt}
\end{eqnarray*}
or, equivalently,
\be
    \bfLf_{(q)} =   T_{(q)} \bfLs_{(q)} \left(T_{(q)}\right)^{-1}
                  + \left(T_{(q)}\right)^{-1} \frac{dT_{(q)}}{dt} .
\Label{Lf-vs-L}
\ee
Next,
we estimate the blocks of
the matrices in
Eq.~(\ref{Lf-vs-L}).
The estimate of $T_{(q)}$
is given in the induction
hypothesis~(\ref{T(q)-exp});
its inverse satisfies a similar estimate,
\begin{eqnarray}
    \left(T_{(q)}\right)^{-1}
    (\,\cdot\,,\psif_{(q)},\eps)
    =
    \left(\begin{array}{cc}
            I_n + {\cal O}(\eps^2)
          & - \eps^q T^{12}_{(q,q)}
          + {\cal O}(\eps^{q+1})\\
          - \eps T^{21}_{(q,1)}
          + {\cal O}(\eps^2)
          & I_m + {\cal O}(\eps^2)
         \end{array}
    \right) .
\Label{Tinv-exp}
\end{eqnarray}
Also,
the induction hypothesis~(\ref{psif(q)-exp})
and
Theorem~\ref{1step-Main-thm}
guarantee that
$\psif_{(q)}
= \psis_{(q)}
+ \mathcal{O}(\eps^{q+1})$,
so
the expansions of
$\bfLs_{(q)}(y,\psis_{(q)},\eps)$
and
$\bfLs_{(q)}(y,\psif_{(q)},\eps)$
are equal up to and including
terms of $\mathcal{O}(\eps^{q})$.
It follows from
Eq.~(\ref{L-order})
that
\be
    \bfLs_{(q)}(\,\cdot\,,\psif_{(q)},\eps)
    =
    \left(\begin{array}{cc}
            \bfLf^{11}_{(0,0)} + {\cal O}(\eps)
          & \eps^q \bfLs^{12}_{(q,q)}
          + {\cal O}(\eps^{q+1})\\
            \eps \bfLf^{21}_{(0,1)}
          + {\cal O}(\eps^2)
          & \eps\bfLs^{22}_{(q,1)} + {\cal O}(\eps^2)
         \end{array}
    \right).
\Label{L(q)-exp}
\ee
Taking $V=T_{(q)}$
in Lemma~\ref{DUg-lem},
we conclude from Eq.~(\ref{T(q)-exp})
that
\be
    DT_{(q)} g
    =
    \left(\begin{array}{cc}
            {\cal O}(\eps^3)
          & {\cal O}(\eps^{q+1})\\
            {\cal O}(\eps^2)
          & {\cal O}(\eps^3)
         \end{array}
    \right) .
\Label{DTg-exp}
\ee
The desired estimate of
$\bfLf_{(q)}$
now follows immediately
from Eqs.~(\ref{T(q)-exp}),
(\ref{Tinv-exp}),
(\ref{L(q)-exp}), and
(\ref{DTg-exp}),
\begin{eqnarray}
    \bfLf^{11}_{(q)}
    &=&
      \bfLf^{11}_{(0,0)}
    + {\cal O}(\eps) ,
\Label{Lf11-exp}\\
    \bfLf^{12}_{(q)}
    &=&
      \eps^q 
      [\bfLs^{12}_{(q,q)}
            - \bfLf^{11}_{(0,0)} T^{12}_{(q,q)} ]
      + {\cal O}(\eps^{q+1}) ,
\Label{Lf12-exp}\\
    \bfLf^{21}_{(q)}
    &=&
    \eps
    [\bfLf^{21}_{(0,1)}
            + T^{21}_{(q,1)} \bfLf^{11}_{(0,0)} ]
    + {\cal O}(\eps^{2}) ,
\Label{Lf21-exp}\\
    \bfLf^{22}_{(q)}
    &=&
    \eps 
    \bfLs^{22}_{(q,1)} 
    + {\cal O}(\eps^2) .
\Label{Lf22-exp}
\end{eqnarray}
The definition~(\ref{U(q),L(q)})
and
Eqs.~(\ref{Lf11-exp})
and
(\ref{Lf12-exp})
imply that
$\Uf_{(q)} = \mathcal{O}(\eps^q)$,
with the leading-order coefficient
given by
\begin{eqnarray}
    \Uf_{(q,q)}
    =
    \left(\bfLf^{11}_{(q,0)}\right)^{-1}
    \bfLf^{12}_{(q,q)}
    =
    \Us_{(q,q)} - T^{12}_{(q,q)} .
\Label{Uf-vs-U}
\end{eqnarray}
Furthermore,
the definition~(\ref{U(q),L(q)})
and
Eqs.~(\ref{Lf11-exp}) 
and 
(\ref{Lf21-exp})
imply that
\begin{eqnarray}
    \Lf_{(q)}
    =
    \bfLf^{21}_{(q)}
    \left(\bfLf^{11}_{(q)}\right)^{-1}
    =
    \mathcal{O}(\eps) .
\Label{Lf-order}
\end{eqnarray}
Finally, we observe that,
to leading order,
the blocks of
$T_{(q)}(\cdot\,,\psif_{(q)},\eps)$
are all equal
to the corresponding blocks of
$T_{(q)}(\cdot\,,\psif_{(q-1)},\eps)$.
The latter are given by
the induction hypothesis~(\ref{T(q)-exp}).

We are now ready to
estimate the size of
the blocks of
$T_{(q+1)}(\cdot\,,\psif_{(q)},\eps)$.

The update formula~(\ref{T-rec}) gives
$T^{11}_{(q+1)}
=
T^{11}_{(q)}
+ \Uf_{(q)} T^{21}_{(q)}$.
According to the induction hypothesis,
$T^{11}_{(q)}=I_n +\mathcal{O}(\eps^2)$
and
$T^{21}_{(q)}=\mathcal{O}(\eps)$.
Furthermore,
$\Uf_{(q)}=\mathcal{O}(\eps^q)$,
by Eq.~(\ref{Uf-vs-U}).
Thus,
$T^{11}_{(q+1)}
=I_n+\mathcal{O}(\eps^2)$,
as desired.

The update formula~(\ref{T-rec}) also gives
$T^{12}_{(q+1)}
=
T^{12}_{(q)}
- T^{11}_{(q)} \Us_{(q)}
+ \Uf_{(q)} T^{22}_{(q)}
- \Uf_{(q)} T^{21}_{(q)} \Us_{(q)}$.
According to the induction hypothesis,
$T^{12}_{(q)}
= \mathcal{O}(\eps^q)$,
$T^{11}_{(q)}
= I_n + \mathcal{O}(\eps^2)$,
$T^{21}_{(q)}
= \mathcal{O}(\eps)$,
and
$T^{22}_{(q)}
= I_m + \mathcal{O}(\eps^2)$.
Furthermore,
$\Uf_{(q)}=\mathcal{O}(\eps^q)$,
by Eq.~(\ref{Uf-vs-U}),
and
$\Us_{(q)}=\mathcal{O}(\eps^q)$,
by Eq.~(\ref{U(q-1)-order}).
Thus, the terms in
the formula for $T^{12}_{(q+1)}$
are all at least
$\mathcal{O}(\eps^q)$.
The same is then true for $T^{12}_{(q+1)}$.
We will now show that
$T^{12}_{(q+1)}$ is, in fact,
at least $\mathcal{O}(\eps^{q+1})$
by showing that 
$T^{12}_{(q+1,q)}=0$.
To leading order,
the update formula for $ T^{12}_{(q+1)}$ is
\be
    T^{12}_{(q+1,q)} = T^{12}_{(q,q)}
                - \Us_{(q,q)}
                + \Uf_{(q,q)}.
\Label{T-exp-lem}
\ee
Equation~(\ref{Uf-vs-U})
implies that
the right member of
(\ref{T-exp-lem})
vanishes.
Therefore,
$T_{(q+1,q)}=0$,
as desired.
We emphasize again that
the choice of $U_{(q)}$
is central to the working of
the CSP method.

Next, the update formula~(\ref{T-rec}) gives
$T^{21}_{(q+1)}
=
T^{21}_{(q)}
- \Lf_{(q)} \Uf_{(q)} T^{21}_{(q)}
- \Lf_{(q)} T^{11}_{(q)}$.
According to the induction hypothesis,
$T^{21}_{(q)}
=\mathcal{O}(\eps)$
and
$T^{11}_{(q)}
=I_n + \mathcal{O}(\eps^2)$.
Furthermore,
$\Uf_{(q)}=\mathcal{O}(\eps^q)$ and
$\Lf_{(q)}=\mathcal{O}(\eps)$,
by Eqs.~(\ref{Uf-vs-U})
and~(\ref{Lf-order}).
Thus, the terms in
the update formula for $T^{21}_{(q+1)}$
are all at least
$\mathcal{O}(\eps)$.
Hence,
$T^{21}_{(q)}$
is also at least
$\mathcal{O}(\eps)$,
as desired.

Lastly,
the update formula~(\ref{T-rec}) gives
$T^{22}_{(q+1)}   
=      
T^{22}_{(q)}
- \Lf_{(q)} T^{12}_{(q)}
- T^{21}_{(q)} \Us_{(q)}
+ \Lf_{(q)} T^{11}_{(q)} \Us_{(q)}
- \Lf_{(q)} \Uf_{(q)} T^{22}_{(q)}
+  \Lf_{(q)} \Uf_{(q)} T^{21}_{(q)} \Us_{(q)}$.
According to the induction hypothesis,
$T^{22}_{(q)}
=I_m + \mathcal{O}(\eps^2)$.
The remaining terms
have already been shown
to be at least
$\mathcal{O}(\eps^2)$.
Hence,
$T^{22}_{(q+1)}
=I_m+\mathcal{O}(\eps^2)$.

The proof of
Theorem~\ref{t-CSPMq}
is complete.
%

%%%%%%%%%%%%%%%%%%%%%%%%%%%%%%%%%%%%%%%%%%%%%%%%%%%%%%%%%%%%%%%%%%%%%%%
%%%%%%%%%%%%%%%%%%%%%%%%%%%%%%%%%%%%%%%%%%%%%%%%%%%%%%%%%%%%%%%%%%%%%%%
%%%%%%%%%%%%%%%%%%%%%%%%%%%%%%%%%%%%%%%%%%%%%%%%%%%%%%%%%%%%%%%%%%%%%%%

\section{The Michaelis--Menten--Henri Reaction \label{s-MMH}}
\setcounter{equation}{0}
In this section, we apply the CSP method
to the Michaelis--Menten--Henri (MMH) mechanism
of enzyme kinetics to illustrate
Theorems~\ref{t-CSPMq} and~\ref{1step-Main-thm}.
We consider
the planar system of ODEs
for a slow variable $s$ and
a fast variable $c$,
\begin{eqnarray}
    s' &=& \eps (- s + (s + \kappa - \lambda)c) , \Label{MMH-s} \\
    c' &=& s - (s + \kappa)c . \Label{MMH-c}
\end{eqnarray}
The parameters 
satisfy the inequalities
$0 < \eps \ll 1$ and $\kappa > \lambda > 0$.
Only nonnegative values of $s$ and $c$ are relevant.
The system of Eqs.~(\ref{MMH-s})--(\ref{MMH-c})
is of the form~(\ref{eq-y})--(\ref{eq-z})
with $m=1$, $n=1$,
$y=s$, $z=c$,
$g_1 = -s+(s+\kappa-\lambda)c$,
and
$g_2 = s-(s+\kappa)c$.

In the limit as $\eps \downarrow 0$,
the dynamics of the MMH equations
are confined to the reduced slow manifold
\be
    \Mo = \{ (c, s) : c = \frac{s}{s+\kappa} , s \ge 0 \} . \quad
    \Label{MMH-Mo}
\ee
The manifold $\Mo$ is normally hyperbolic,
so according to Theorem~\ref{t-Fenichel}
there exists,
for all sufficiently small $\eps$,
a slow manifold $\Me$ 
that is $\mathcal{O}(\eps)$ close
to $\Mo$ on any compact set. 
Moreover, $\Me$ is the graph
of a function $h_\eps$,
\be
    \Me = \{ (c, s) : c = h_\eps (s) , s \ge 0 \} ,
    \Label{MMH-Me}
\ee
and $h_\eps$ admits an asymptotic expansion 
$h_\eps = h_0 + \eps h_1 + \eps^2 h_2 + \cdots \,$.
The coefficients are found from
the invariance equation,
\be
  s - (s + \kappa) h_\eps (s)
  =
  \eps h_\eps' (s) (-s + (s + \kappa - \lambda) h_\eps (s)) .
\ee 
The first few coefficients are
\be
    h_0 (s) = \frac{s}{s+\kappa} , \quad
    h_1(s) = \frac{\kappa\lambda s}{(s+\kappa)^4} , \quad
    h_2(s) = \frac{\kappa\lambda s(2\kappa\lambda-3\lambda s-\kappa s-\kappa^2)}{(s+\kappa)^7} .
    \Label{Me-MMH}
\ee

\subsection{Application of the One-Step CSP Method
\label{ss-MMH}}
Both the one-step and two-step CSP methods
start from the same initial basis.
We choose the stoichiometric vectors
as the basis vectors, so
\be
    \Af^{(0)} 
    = 
    (\Af_1^{(0)}, \, \Af_2^{(0)})
    = 
    \left(
          \begin{array}{cc}
          0 & 1 \\
          1 & 0
          \end{array}
    \right) , \quad
    \Bf_{(0)}
    =
    \left(
          \begin{array}{cc}
          \Bf^1_{(0)} \\
          \Bf^2_{(0)}
          \end{array}
     \right)
     =
     \left(
          \begin{array}{cc}
          0 & 1 \\
          1 & 0
          \end{array}
      \right) .
\Label{1/A,B}
\ee
The CSP condition
$B^1_{(0)} g = 0$
is satisfied if
$c = h_0 (s)$,
so the CSP manifolds
$\Ks^{(0)}$ and $\Kf^{(0)}$
coincide with $\Mo$.
With this choice of
initial basis, we have
\be
    \bfLf_{(0)}
    = 
    \Bf_{(0)} (Dg) \Af^{(0)} 
    = 
    \left(
          \begin{array}{cc}
            - (s+\kappa) 
          & - (c-1) \\
            \eps(s+\kappa-\lambda) 
          & \eps(c-1)
          \end{array}
    \right) .
\Label{1/MMH-Lambda(0)}
\ee
\paragraph{First iteration.}
At any point $(s, c)$,
we have
\be
    \As^{(1)}
=
    \left(
          \begin{array}{cc}
          0 & 1\\
          1 & - \frac{c-1}{s+\kappa}
          \end{array}
    \right) , \qquad
    \Bs_{(1)}
=
    \left(
          \begin{array}{cc}
          \frac{c-1}{s+\kappa} & 1\\
          1 & 0
          \end{array}
    \right) .
\Label{1/a1, b1 genl pt}
\ee
On $\Ks^{(0)}$, these expressions reduce to
\be
    \As^{(1)}
=
    \left(
          \begin{array}{cc}
          0 & 1\\
          1 & \frac{\kappa}{(s+\kappa)^2}
          \end{array}
    \right) , \qquad
    \Bs_{(1)}
=
    \left(
          \begin{array}{cc}
          \frac{-\kappa}{(s+\kappa)^2} & 1\\
          1 & 0
          \end{array}
    \right) .
\Label{1/a1, b1}
\ee
The CSP condition,
\be
    \Bs^1_{(1)} g
=  
      s
    - (s+\kappa)c
    - \eps\frac{\kappa(-s+(s+\kappa-\lambda)c)}{(s+\kappa)^2}
=
    0 ,
\Label{1/f11}
\ee
is satisfied if
\be
    c 
%=
%   \frac{  s
%       + \eps\frac{\kappa s}{(s+\kappa)^2}
%      }
%   {  (s+\kappa)
%    + \eps\frac{\kappa(s+\kappa-\lambda)}{(s+\kappa)^2}
%   }
= 
      \frac{s}{s+\kappa}
    + \eps\frac{\kappa\lambda s}
        {(s+\kappa)^4}
    - \eps^2\frac{\kappa^2\lambda s(s+\kappa-\lambda)}
          {(s+\kappa)^7}
    + \mathcal{O}(\eps^3) . 
\Label{1/proj2}
\ee
Comparing this result
with Eq.~(\ref{Me-MMH}),
we see that
the asymptotic expansions
of $\Ks^{(1)}$ and $\Me$
coincide
up to and including $\mathcal{O}(\eps)$ terms,
in accordance with
Theorem~\ref{1step-Main-thm}
for $q=1$;
however, the $\mathcal{O}(\eps^2)$ terms
differ at this stage.

\paragraph{Second iteration.}
The blocks of
$\bfLs_{(1)}$
are
\begin{eqnarray}
\hspace{-1em}
    \bfLs^{11}_{(1)}
&=&
    - (s+\kappa)
    + \eps \frac{(s+\kappa-\lambda)(c-1)}{s+\kappa} ,
\Label{1/L11 genl pt} \\
    \bfLs^{12}_{(1)}
&=&
      \frac{s}{s+\kappa}-c
    + \eps 
      \frac{(c-1)
            [\lambda(c-1)
          - (-s+(s+\kappa-\lambda)c)]}
      {(s+\kappa)^2} ,
\Label{1/L12 genl pt} \\
    \bfLs^{21}_{(1)}
&=&
    \eps (s+\kappa-\lambda) , \qquad
    \bfLs^{22}_{(1)}
=
    \eps \frac{\lambda (c-1)}{s+\kappa} .
\Label{1/L21,22 genl pt}
\end{eqnarray}
On $\Ks^{(1)}$,
the blocks reduce to
\begin{eqnarray}
    \bfLs^{11}_{(1)}
&=&
    - (s+\kappa)
    - \eps \frac{\kappa (s+\kappa-\lambda)}{(s+\kappa)^2}
    + \eps^2 \frac{\kappa\lambda s(s+\kappa-\lambda)}{(s+\kappa)^5} ,
\Label{1/L11} \\
    \bfLs^{12}_{(1)}
&=&
      \eps \frac{\kappa\lambda(\kappa-2s)}{(s+\kappa)^4}
    + \eps^2 
      \frac{\kappa\lambda s(2\kappa(s+\kappa-2\lambda)+\lambda s)}
      {(s+\kappa)^7} ,
\Label{1/L12} \\
    \bfLs^{21}_{(1)}
&=&
    \eps (s+\kappa-\lambda) , \qquad
    \bfLs^{22}_{(1)}
=
    - \eps \frac{\kappa\lambda}{(s+\kappa)^2}
    + \eps^2 \frac{\kappa\lambda^2 s}{(s+\kappa)^5} .
\Label{1/L21,22}
\end{eqnarray}
The second update is
\begin{eqnarray}
    \As_1^{(2)} 
&=&
    \left(
          \begin{array}{c}
            0\\
            1
          \end{array}
    \right) ,
\Label{1/a12}\\
    \As_2^{(2)} 
&=&
    \left(
          \begin{array}{c}
            1\\
            \frac{\kappa}{(s+\kappa)^2}
          \end{array}
    \right)
    +
    \eps
    \left(
         \begin{array}{cc}
           0\\
           \frac{\kappa\lambda(\kappa-3s)}{(s+\kappa)^5}
         \end{array}
    \right)
\nonumber\\
    &&\mbox{}+
    \eps^2
    \left(
          \begin{array}{cc}
            0\\
            \frac{\kappa\lambda[\kappa(5s-\kappa)(s+\kappa-\lambda)+\lambda s(s-2\kappa)]}
            {(s+\kappa)^8}
          \end{array}
    \right)
    + \mathcal{O}(\eps^3) ,
\Label{1/a22}\\
%
%    \As^{(2)} 
%&=&
%    \left(
%          \begin{array}{cc}
%            0 & 1\\
%            1 & \frac{\kappa}{(s+\kappa)^2}
%          \end{array}
%    \right)
%    +
%    \eps
%    \left(
%         \begin{array}{cc}
%           0 & 0\\
%           0 & \frac{\kappa\lambda(\kappa-3s)}{(s+\kappa)^5}
%         \end{array}
%    \right)
%\nonumber\\
%    &&\mbox{}+
%    \eps^2
%    \left(
%          \begin{array}{cc}
%            0 & 0\\
%            0 & \frac{\kappa\lambda(\kappa(5s-\kappa)(s+\kappa-\lambda)+\lambda s(s-2\kappa))}
%            {(s+\kappa)^8}
%          \end{array}
%    \right)
%    +
%    \mathcal{O}(\eps^3) ,
%\Label{1/a2}\\
%
    \Bs^1_{(2)} 
&=&
    \left(
          - \frac{\kappa}{(s+\kappa)^2},\,1
    \right)
    +
    \eps
    \left(
            - \frac{\kappa\lambda(\kappa-3s)}
              {(s+\kappa)^5},\, 0
    \right)
    \nonumber\\
    &&\mbox{}
    + \eps^2
    \left(
            - \frac{\kappa\lambda[\kappa(5s-\kappa)(s+\kappa-\lambda)+\lambda s(s-2\kappa)]}
              {(s+\kappa)^8},\, 0
    \right)
    \nonumber\\
    &&\mbox{}
    + \mathcal{O}(\eps^3) ,
\Label{1/b12}\\
    \Bs^2_{(2)} 
&=&
    \left(
            1,\, 0
    \right) .
\Label{1/b22}
%
%    \Bs_{(2)} 
%&=&
%    \left(
%          \begin{array}{cc}
%          - \frac{\kappa}{(s+\kappa)^2} & 1\\
%            1 & 0
%          \end{array}
%    \right)
%    +
%    \eps
%    \left(
%          \begin{array}{cc}
%            - \frac{\kappa\lambda(\kappa-3s)}
%              {(s+\kappa)^5} & 0\\
%            0 & 0
%          \end{array}
%    \right)
%    \nonumber\\
%    &&\mbox{}
%    + \eps^2
%    \left(
%          \begin{array}{cc}
%            - \frac{\kappa\lambda(\kappa(5s-\kappa)(s+\kappa-\lambda)+\lambda s(s-2\kappa))}
%              {(s+\kappa)^8} & 0\\
%            0 & 0
%         \end{array}
%    \right)
%    + \mathcal{O}(\eps^3) .
%\Label{1/b2}
\end{eqnarray}
The CSP condition
\begin{eqnarray}
    \Bs^1_{(2)} g
&=&
      s
    - (s+\kappa)c
    - \eps\frac{\kappa(-s+(s+\kappa-\lambda)c)}{(s+\kappa)^2}
      \nonumber\\
    &&\mbox{}
    + \eps^2\kappa\lambda
      \frac{(3s-\kappa)(-s+(s+\kappa-\lambda)c)}{(s+\kappa)^5}
    + \mathcal{O}(\eps^3)  \nonumber \\
&=& 0 ,
\Label{1/f12}
\end{eqnarray}
is satisfied if
\begin{eqnarray}
    c 
%&=& 
%   \frac{
%   s
%
%   \eps\frac{\kappa s}{(s+\kappa)^2}
%-
%   \eps^2\kappa\lambda
%   \frac{(3s-\kappa)s}
%   {(s+\kappa)^5}
%+
%   \mathcal{O}(\eps^3)
%   }
%   {
%   (s+\kappa)
%+
%   \eps\frac{\kappa(s+\kappa-\lambda)}{(s+\kappa)^2}
%-
%    \eps^2\kappa\lambda
%    \frac{(s+\kappa-\lambda)(3s-\kappa)}
%    {(s+\kappa)^5}
%+
%   \mathcal{O}(\eps^3)
%   }
%\nonumber \\
&=& 
    \frac{s}{s+\kappa}
+
    \eps\frac{\kappa\lambda s}{(s+\kappa)^4}
+    
    \eps^2
     \frac{\kappa\lambda s(2\kappa\lambda-3\lambda s-\kappa s-\kappa^2)}{(s+\kappa)^7}+\mathcal{O}(\eps^3) . 
\Label{1/proj3}
\end{eqnarray}
Comparing this result with
Eq.~(\ref{Me-MMH}),
we see that
the asymptotic expansions
of $\Ks^{(2)}$ and $\Me$
coincide
up to and including $\mathcal{O}(\eps^2)$ terms,
in accordance with
Theorem~\ref{1step-Main-thm}
for $q=2$.

\subsection{Application of the Full CSP Method
\label{ss-MMH-full}}

\paragraph{First iteration.}
At any point $(s,c)$,
we have
\be
    \Af_1^{(1)}
=
    \left(
          \begin{array}{c}
          0\\
          1
          \end{array}
    \right)
    -
    \eps\frac{s+\kappa-\lambda}{s+\kappa}
    \left(
          \begin{array}{c}
           1\\
          \mbox{}- \frac{c-1}{s+\kappa}
          \end{array}
    \right) , \qquad
    \Af_2^{(1)} 
=
    \left(
          \begin{array}{c}
            1 \\
          \mbox{}- \frac{c-1}{s+\kappa}
          \end{array}
    \right) ,
\Label{a1 genl pt}
\ee
\be
    \Bf^1_{(1)} 
=
    \left(
          \frac{c-1}{s+\kappa}, \, 1
    \right), \qquad
    \Bf^2_{(1)} 
=
    \left(
          1 , \, 0
    \right)
    +
    \eps\frac{s+\kappa-\lambda}{s+\kappa}
    \left(
          \frac{c-1}{s+\kappa}, \, 1
    \right) .
\Label{b1 genl pt}
\ee
On $\Kf^{(0)}$,
these quantities reduce to
\be
    \Af_1^{(1)} 
=
    \left(
          \begin{array}{c}
          0\\
          1
          \end{array}
    \right)
    -
    \eps\frac{s+\kappa-\lambda}{s+\kappa}
    \left(
          \begin{array}{c}
          1\\
          \frac{\kappa}{(s+\kappa)^2}
          \end{array}
    \right), \qquad
    \Af_2^{(1)} 
=
    \left(
          \begin{array}{c}
          1\\
          \frac{\kappa}{(s+\kappa)^2}
          \end{array}
    \right),
\Label{a1}
\ee
\be
    \Bf^1_{(1)} 
=
    \left(
          \mbox{}-\frac{\kappa}{(s+\kappa)^2}, \, 1
    \right) ,  \quad
    \Bf^2_{(1)} 
=
    \left(
          1, \, 0
    \right)
    +
    \eps\frac{s+\kappa-\lambda}{s+\kappa}
    \left(
          \mbox{}- \frac{\kappa}{(s+\kappa)^2}, \, 1
    \right).
\Label{b1}
\ee
The matrix relating $\Bf_{(1)}$
to its one-step counterpart~$\Bs_{(1)}$
is
\be
    T_{(1)}
    =
    \left(
          \begin{array}{cc}
            1 & 0\\
            \eps\frac{s+\kappa-\lambda}{s+\kappa}
          & 1
          \end{array}
    \right),
\Label{T(1)}
\ee
so
$T_{(1)}$
is indeed of the form
(\ref{T(q)-exp})
on $\Kf^{(0)}$.

Equations~(\ref{1/a1, b1})
and
(\ref{b1})
imply that
$\Bf_{(1)}=\Bs_{(1)}$,
so
the CSP condition
yields
$\psif_{(1)}=\psis_{(1)}$.
Thus, after one iteration,
the full CSP method 
also finds
the expansion of ${\cal M}_\eps$ 
up to and including 
$\mathcal{O}(\eps)$ terms.

\paragraph{Second iteration.}
The blocks of 
$\bfLf_{(1)}$
are
\begin{eqnarray}
    \bfLf^{11}_{(1)}
&=&
    - (s+\kappa)
    + \eps \frac{(s+\kappa-\lambda)}{s+\kappa}
      \left[(c-1)
      + (c-\frac{s}{s+\kappa})\right]
    \nonumber\\
    &&\mbox{}
    + \eps^2 \frac{(c-1)(s+\kappa-\lambda)}{(s+\kappa)^3}
      \left[- \lambda (c-1)
       + (-s+(s+\kappa-\lambda)c)\right] ,~~
%    + \mathcal{O}(\eps^3) ,
\Label{L11 genl pt}\\
    \bfLf^{12}_{(1)}
&=&
      \frac{s}{s+\kappa}-c
    + \eps \frac{c-1}{(s+\kappa)^2}
    \left[\lambda (c-1)
     - (-s+(s+\kappa-\lambda)c)\right] ,
%    + \mathcal{O}(\eps^3) , 
\Label{L12 genl pt}\\
    \bfLf^{21}_{(1)}
&=&
    \left.\frac{\eps^2}{(s+\kappa)^2} 
    \right[(c-1)(s+\kappa-\lambda)(s+\kappa-2\lambda)
    \nonumber\\
    &&\mbox{}\left.
            + \lambda (-s+(s+\kappa-\lambda)c)
            + (s+\kappa-\lambda)^2
              \left(
                    c-\frac{s}{s+\kappa}
              \right)\right] ,
%    + \mathcal{O}(\eps^3) ,
\Label{L21 genl pt}\\
    \bfLf^{22}_{(1)}
&=&
    \frac{\eps}{s+\kappa} 
    \left[\lambda (c-1)
          + (s+\kappa-\lambda)(\frac{s}{s+\kappa}-c)\right]
    \nonumber\\
    &&\mbox{}
    + \eps^2 \frac{(c-1)(s+\kappa-\lambda)}{(s+\kappa)^3}
    \left[\lambda (c-1) - (-s+(s+\kappa-\lambda)c)\right] ,
%    + \mathcal{O}(\eps^3) .
\Label{L22 genl pt}
\end{eqnarray}
with remainders of $\mathcal{O}(\eps^3)$.
On $\Kf^{(1)}$,
the blocks reduce to
\begin{eqnarray}
    \bfLf^{11}_{(1)}
&=&
    - (s+\kappa)
    - \eps \frac{\kappa (s+\kappa-\lambda)}{(s+\kappa)^2}
    + \eps^2 \frac{\kappa\lambda(s+\kappa-\lambda)(3s-\kappa)}{(s+\kappa)^5} ,
\Label{L11}\\
    \bfLf^{12}_{(1)}
&=&
      \eps \frac{\kappa\lambda(\kappa-2s)}{(s+\kappa)^4}
    + \eps^2 \frac{\kappa\lambda s(2\kappa(s+\kappa-2\lambda)+\lambda s)}{(s+\kappa)^7} ,
\Label{L12}\\
    \bfLf^{21}_{(1)}
&=&
    - \eps^2 \frac{\kappa(s+\kappa-\lambda)(s+\kappa-2\lambda)+\lambda^2 s}{(s+\kappa)^3} ,
\Label{L21}\\
    \bfLf^{22}_{(1)}
&=&
    - \eps \frac{\kappa\lambda}{(s+\kappa)^2}
    - \eps^2 \frac{\kappa\lambda((2s-\kappa)(s+\kappa-\lambda)-\lambda s)}{(s+\kappa)^5} ,
\Label{L22}
%
%
%
%
%
%     \bfLf_{(1)} 
%&=& 
%      \left(
%            \begin{array}{cc}
%            - (s+\kappa) & 
%               0\\
%               0 & 0
%            \end{array}
%      \right)
%      +
%      \eps\left(
%                \begin{array}{cc}
%                - \kappa\frac{s+\kappa-\lambda}{(s+\kappa)^2} & 
%                  \frac{\kappa\lambda(\kappa-2s)}{(s+\kappa)^4}\\
%                  0 & 
%                - \frac{\kappa\lambda}{(s+\kappa)^2}\\
%                \end{array}
%          \right)
%      \nonumber\\
%      &&\mbox{}
%      + \eps^2\left(
%                    \begin{array}{cc}
%                      \frac{\kappa\lambda(s+\kappa-\lambda)(3s-\kappa)}{(s+\kappa)^5} &
%                      \frac{\kappa\lambda s(2\kappa(s+\kappa-2\lambda)+\lambda s)}{(s+\kappa)^7} \\
%                    - \frac{\kappa(s+\kappa-2\lambda)(s+\kappa-\lambda)+\lambda^2 s}{(s+\kappa)^3} &
%                    - \frac{\kappa\lambda((2s-\kappa)(s+\kappa-\lambda)-\lambda s)}{(s+\kappa)^5}
%                    \end{array}
%              \right)
%      \nonumber\\
%      &&\mbox{}
%      + \mathcal{O}(\eps^3).
%\Label{L11}
\end{eqnarray}
with errors of $\mathcal{O}(\eps^3)$.
The result of the second iteration is
\begin{eqnarray}
    \Af_{11}^{(2)} 
&=&
    - \eps \frac{s+\kappa-\lambda}{s+\kappa}
    + \eps^2 
    \frac{\kappa(s+\kappa-2\lambda)(s+\kappa-\lambda)+\lambda^2 s}{(s+\kappa)^4} ,
\Label{a11}\\
    \Af_{12}^{(2)} 
&=&
      1
    + \eps^2 
    \frac{\kappa\lambda (2s-\kappa)(s+\kappa-\lambda)}{(s+\kappa)^6} ,
\Label{a12}\\
    \Af_{21}^{(2)} 
&=&
      1
    - \eps \frac{\kappa(s+\kappa-\lambda)}{(s+\kappa)^3} 
    \nonumber\\
    &&\mbox{}
    + \eps^2
    \frac{(s+\kappa-\lambda)(\kappa^2(s+\kappa-2\lambda)+\kappa\lambda s)+\kappa\lambda^2 s}
    {(s+\kappa)^6},
\Label{a21}\\
    \Af_{22}^{(2)} 
&=&
      \frac{\kappa}{(s+\kappa)^2}
    + \eps \frac{\kappa\lambda(\kappa-3s)}{(s+\kappa)^5}
    \nonumber\\
    &&\mbox{}
    + \eps^2
    \frac{\kappa^2\lambda (7s-2\kappa)(s+\kappa-\lambda)+\kappa\lambda^2 s (s-2\kappa)}
    {(s+\kappa)^8} ,
\Label{a22}\\
%
%
%
%
%
%    \Af^{(2)} 
%&=&
%    \left(
%          \begin{array}{cc}
%            0 
%          & 1\\
%            1
%          & \frac{\kappa}{(s+\kappa)^2}
%          \end{array}
%    \right)
%    +
%    \eps
%    \left(
%          \begin{array}{cc}
%            - \frac{s+\kappa-\lambda}{s+\kappa}
%          & 0\\
%            - \frac{\kappa(s+\kappa-\lambda)}{(s+\kappa)^3}          
%          & \frac{\kappa\lambda(\kappa-3s)}{(s+\kappa)^5}
%          \end{array}
%    \right)
%    \nonumber\\
%    &&\mbox{}
%    + 
%    \eps^2
%    \left(
%          \begin{array}{cc}
%\frac{\kappa(s+\kappa-2\lambda)(s+\kappa-\lambda)+\lambda^2 s}{(s+\kappa)^4}
%& A_{12}^{(2,2)}\\
%\frac{(s+\kappa-\lambda)
%(\kappa^2(s+\kappa-2\lambda)+\kappa\lambda s)
%+\kappa\lambda^2 s}
%{(s+\kappa)^6}
%& A_{22}^{(2,2)}
%          \end{array}
%    \right)
%    +
%    \mathcal{O}(\eps^3),
%\Label{a2}\\
%
    \Bf^{11}_{(2)} 
&=&
      \frac{-\kappa}{(s+\kappa)^2}
    - \eps \frac{\kappa\lambda(\kappa-3s)}{(s+\kappa)^5}
    \nonumber\\
    &&\mbox{}
    - \eps^2
    \frac{\kappa^2\lambda (7s-2\kappa)(s+\kappa-\lambda)+\kappa\lambda^2 s (s-2\kappa)}
    {(s+\kappa)^8} ,
\Label{b11}\\
    \Bf^{12}_{(2)} 
&=&
      1
    + \eps^2 
    \frac{\kappa\lambda (2s-\kappa)(s+\kappa-\lambda)}{(s+\kappa)^6} ,
\Label{b12}\\
    \Bf^{21}_{(2)} 
&=&
      1
    - \eps \frac{\kappa(s+\kappa-\lambda)}{(s+\kappa)^3}
    \nonumber\\
    &&\mbox{}
    + \eps^2
    \frac{(s+\kappa-\lambda)(\kappa^2(s+\kappa-2\lambda)+\kappa\lambda s)+\kappa\lambda^2 s}
    {(s+\kappa)^6} ,
\Label{b21}\\
    \Bf^{22}_{(2)} 
&=&
      \eps \frac{s+\kappa-\lambda}{s+\kappa}
    - \eps^2 
    \frac{\kappa(s+\kappa-\lambda)(s+\kappa-2\lambda)+\lambda^2 s}
    {(s+\kappa)^4} ,
\Label{b22}
%
%
%
%
%
%    \Bf_{(2)} 
%&=&
%    \left(
%          \begin{array}{cc}
%            \frac{-\kappa}{(s+\kappa)^2} 
%          & 1\\
%            1
%          & 0
%          \end{array}
%    \right)
%%
%    +
%%
%    \eps
%    \left(
%          \begin{array}{cc}
%            - \frac{\kappa\lambda(\kappa-3s)}{(s+\kappa)^5}
%          & 0\\
%            - \frac{\kappa(s+\kappa-\lambda)}
%              {(s+\kappa)^3}
%          & \frac{s+\kappa-\lambda}{(s+\kappa)^3}
%          \end{array}    
%    \right)
%    \nonumber\\
%    &&\mbox{}
%    + 
%    \eps^2
%    \left(
%          \begin{array}{cc}
%            B^{11}_{(2,2)}
%          & B^{12}_{(2,2)}\\
%            \frac{(s+\kappa-\lambda)
%            (\kappa^2(s+\kappa-2\lambda)+\kappa\lambda s)
%            + \kappa\lambda^2 s}
%            {(s+\kappa)^6}
%          & - \frac{\kappa(s+\kappa-\lambda)(s+\kappa-2\lambda)+\lambda^2 s}
%              {(s+\kappa)^4}
%          \end{array}
%    \right)
%    \nonumber\\
%    &&\mbox{}
%    + \mathcal{O}(\eps^3) .
%\Label{b2}
\end{eqnarray}
up to and including terms of $\mathcal{O}(\eps^2)$.
Also,
on $\Kf^{(1)}$,
\begin{eqnarray*}
    T^{11}_{(2)}
&=&
      1 
    + \eps^2
    \frac{\kappa\lambda(s+\kappa-\lambda)(2s-\kappa)}{(s+\kappa)^6} ,\\
    T^{12}_{(2)}
&=&
    0 ,\\
    T^{21}_{(2)}
&=&
      \eps \frac{(s+\kappa-\lambda)}{s+\kappa}
    - \eps^2
    \frac{\kappa(s+\kappa-\lambda)(s+\kappa-2\lambda)+\lambda^2 s}
    {(s+\kappa)^4} ,\\
    T^{22}_{(2)}
&=&
      1 
    - \eps^2
    \frac{\kappa\lambda(2s-\kappa)(s+\kappa-\lambda)}{(s+\kappa)^8} ,
\end{eqnarray*}
with remainders of $\mathcal{O}(\eps^3)$.
Thus, $T_{(2)}$ is indeed of the form~(\ref{T(q)-exp}) on $\Kf^{(1)}$.

The CSP condition
\begin{eqnarray}
    \Bf^1_{(2)} g
&=&
      s
    - (s+\kappa)c
    - \eps\frac{\kappa(-s+(s+\kappa-\lambda)c)}{(s+\kappa)^2}
    \nonumber\\
    &&\mbox{}
    + \eps^2\kappa\lambda
    \left(  \frac{(3s-\kappa)(-s+(s+\kappa-\lambda)c)}
            {(s+\kappa)^5}
    \right.
    \nonumber\\
    &&\left.\hspace{4em}\mbox{}
          + \frac{(2s-\kappa)(s+\kappa-\lambda)(s-(s+\kappa)c)}
            {(s+\kappa)^6}
    \right)
    \nonumber\\
    &&\mbox{}
    + \mathcal{O}(\eps^3) \nonumber \\
&=& 0 ,
\Label{f12}
\end{eqnarray}
is satisfied if
\begin{eqnarray}
    c 
%&=& 
%    \frac{s+\eps\frac{\kappa s}{(s+\kappa)^2}
%+
%   \eps^2\kappa\lambda s\frac{\lambda(\kappa-2s)-(s+\kappa)s}{(s+\kappa)^6}
%+
%    \mathcal{O}(\eps^3)}
%%
%   {(s+\kappa)+\eps\frac{\kappa(s+\kappa-\lambda)}{(s+\kappa)^2}
%-
%   \eps^2\kappa\lambda s\frac{(s+\kappa-\lambda)(s+\kappa)}{(s+\kappa)^6}
%+
%   \mathcal{O}(\eps^3)} \nonumber \\
&=& 
    \frac{s}{s+\kappa}
+
    \eps\frac{\kappa\lambda s}{(s+\kappa)^4}
+    
    \eps^2
     \frac{\kappa\lambda s(2\kappa\lambda-3\lambda s-\kappa s-\kappa^2)}{(s+\kappa)^7}+\mathcal{O}(\eps^3). 
\Label{proj3}
\end{eqnarray}
Therefore, after two iterations,
the full CSP method finds
the expansion of $\Me$ 
up to and including $\mathcal{O}(\eps^2)$ terms.

%\paragraph{Remark.}
%The slow basis vector
%$\Af_2^{(2)}$  
%is tangent to $\Me$
%up to and including
%$\mathcal{O}(\eps)$ terms
%at all points of $\Kf^{(1)}$.

\subsection{The Second Step and the Fast Fibers of $\mathbf{\Me}$}
The preceding analysis of 
the full CSP method
shows that,
at the $q$th iteration, 
the second step
alters only the terms of
${\mathcal O}(\eps^{q+1})$,
leaving the terms of
${\mathcal O} (1)$
through
${\mathcal O}(\eps^q)$
invariant.
Here, we observe that
the second step also plays 
a constructive role
for the dynamics 
near the slow manifold.
As can be seen in the MMH~example,
the second step yields the asymptotic expansions of
the tangent spaces of the fast fibers at their basepoints
up to and including terms of $\mathcal{O}(\eps^{q+1})$,
at least for $q=0, 1$, and $2$.
This additional information is contained in
the columns of $\Af^{(q)}_1$.
We remark here that
this property is not shared by 
the one-step CSP method, 
since the columns of $\As^{(q)}_1$
remain tangent to the fast fibers 
at their basepoints only to leading order
after each iteration.
Details about the fast fibers
and their tangent spaces
will be presented in a future publication.

\section{Relation between CSPM and ILDM \label{s-discussion}}
\setcounter{equation}{0}
The CSP iteration procedure
is designed to diagonalize
the Lie bracket $[\, \cdot\, , g]$.
At each iteration,
the then-current basis is updated
in such a way that
$[\, \cdot\, , g]$
is block-diagonalized
to the next-higher order in~$\eps$.
Thus, each iteration improves
the quality of  the basis of
the orthogonal complement of
the tangent space.
The CSPM is defined as
the locus of points where
the vector field is orthogonal to
that orthogonal complement.

The ILDM method works, instead, with
the Jacobian, $Dg$, of
Eqs.(\ref{eq-y-slow})--(\ref{eq-z-slow}).
A Schur decomposition transforms
$Dg$ into upper
triangular form,
\be
  Dg = Q N Q' , \quad
  N = \left( \begin{array}{cc} N_s & N_{sf} \\ 0 & N_f \end{array} \right) ,
\ee
where $Q = (Q_s \; Q_f)$ is unitary.
The eigenvalues of $Dg$ appear on the diagonal
of $N$ in descending order of their real parts,
from least negative in the upper left
to most negative in the lower right.
The first $m$ Schur vectors
(the columns of $Q_s$)
form an orthogonal basis
of the slow subspace
and the remaining $n$ Schur vectors
(the columns of $Q_f$)
an orthogonal basis of 
the orthogonal complement of
the slow subspace.
The vector field $g$
is entirely in the slow subspace
if it is orthogonal to
this orthogonal complement---that is,
if
\be
  Q_f' g = 0 .
\ee
This equation defines the ILDM;
see~\cite[Section 3]{KK-2002}.

As we showed in~\cite{KK-2002},
the ILDM is only a first-order
approximation to $\Me$.
The error is always $\mathcal{O} (\eps^2)$
unless $\Mo$ is linear.
The error can be traced back
to the choice of the operator.
The tangent space is
a left-invariant subspace
of the Jacobian
only to leading order,
so putting $Dg$ in upper triangular form
yields the orthogonal complement
only to leading order.
Since the linearized system is only
an approximation of the original
ODEs~(\ref{eq-y})--(\ref{eq-z}),
this choice does not
produce an exact result
unless $g$ is linear.
The success of the CSP method in 
approximating the slow manifold is due to
the fact that the ODEs for the amplitudes~$f$
are equivalent to the
ODEs~(\ref{eq-y})--(\ref{eq-z}).
That is, the full nonlinearity is retained.

The time-derivative term
in the definition~(\ref{Lambda-alt})
must be included in the evaluation
of~$\Lambda$;
otherwise, the accuracy of the CSP method
is compromised.
In fact, such an omission results in 
implementing the ILDM 
rather than the CSP method,
which may be seen as follows.
With our initial choice of
a point-independent basis $\Af^{(0)}$,
the matrix $\bfLf_{(0)}$ is similar to $Dg$;
see Eq.~(\ref{Lambda(q)}).
The omission of the term $(dB_{(q)}/dt) A^{(q)}$
in the calculation of $\Lambda_{(q)}$,
for $q=1,2,\ldots\,$,
would lead to the formula
$\bfLf_{(q)} = 
(I + \Ps_{(q)})
\Bf_{(0)} (Dg) \Af^{(0)}
(I - \Ps_{(q)})$,
which would imply that
$\Lambda_{(q)}$
is similar to $Dg$.
Therefore, 
the one-step CSP method
would put $Dg$,
rather than $\Lambda$,
in lower-triangular form,
just like the ILDM method.
After the second iteration,
one would make an error
(proportional to the curvature of $\Mo$)
at ${\mathcal O}(\eps^2)$, which
subsequent iterations would not remove.
The MMH~example in Section~\ref{s-MMH}
illustrates these observations.

\begin{center}
\textbf{ACKNOWLEDGMENTS}
\end{center}
\noindent
We thank
Harvey Lam and Dimitris Goussis
for generously sharing
their insights into
the CSP method
and our colleague
Michael Davis
for stimulating conversations
in the course of this investigation.

\noindent
The work of H.~K.\ was supported by
the Mathematical, Information,
and Computational Sciences Division
subprogram of the Office of
Advanced Scientific Computing Research,
Office of Science,
U.S.~Department of Energy,
under Contract W-31-109-Eng-38.
The work of T.~K. \ and A.~Z. \ was supported in part 
by the Division of Mathematical Sciences
of the National Science Foundation via grant NSF-0072596.

%%%%%%%%%%%%%%%%%%%%%%%%%%%%%%%%%%%%%%%%%%%%%%%%%%%%%%%%%%%%%%%%%%%%%%%
%%%%%%%%%%%%%%%%%%%%%%%%%%%%%%%%%%%%%%%%%%%%%%%%%%%%%%%%%%%%%%%%%%%%%%%
%%%%%%%%%%%%%%%%%%%%%%%%%%%%%%%%%%%%%%%%%%%%%%%%%%%%%%%%%%%%%%%%%%%%%%%

\appendix
\section{Auxiliary Lemmas \label{s-lemmas}}
\begin{LEMMA}
\Label{Dhq-lem}
The quantity $Dh_q$
is given by 
the formula
\begin{eqnarray}
    Dh_q
    &=&\mbox{}
    - ((D_z g_2)_0)^{-1} 
    \left[
    \left(D_y g_2\right)_q 
    + \sum_{i=0}^{q-1} \left(D_z g_2\right)_{q-i} Dh_i
    - \sum_{\ell=0}^{q-1} D^2 h_\ell g_{1,q-\ell-1}
    \right.\nonumber\\
    &&\mbox{} 
    \left.
    - \sum_{\ell=0}^{q-1} 
    (Dh_\ell) (D_y g_1)_{q-1-\ell}
    - \sum_{i=0}^{q-1}\sum_{\ell=0}^{q-1-i} 
    Dh_\ell (D_z g_1)_{q-1-i-\ell} Dh_i
    \right].
\Label{Dhq}
\end{eqnarray}
\end{LEMMA}
\begin{PROOF}
The coefficient
$h_q$
is found from
the $\mathcal{O}(\eps^q)$
terms in
the invariance equation
(\ref{inveq}),
\be
    g_{2,q} = \sum_{\ell=0}^{q-1} (Dh_\ell) g_{1,q-\ell-1} .
\Label{inveq,q}
\ee
Taking the total derivative
with respect to
$y$
of
both sides
of
(\ref{inveq,q}),
we find
\begin{eqnarray}
    \frac{d}{dy}g_{2,q}
    &=&\mbox{}
    \sum_{\ell=0}^{q-1} (D^2 h_\ell) g_{1,q-\ell-1}
    +
    \sum_{\ell=0}^{q-1} (Dh_\ell) \frac{d}{dy} g_{1,q-\ell-1}.
\Label{Dhq-aux1}
\end{eqnarray}
The operations of
taking
the total derivative
with respect to
$y$
and 
expanding
with respect to
$\eps$
commute,
because
the Fenichel theory
guarantees
$C^r$
smoothness in
$\eps$ and $y$
for each
$r$.
Therefore,
\begin{eqnarray}
    \frac{d}{dy}g_{2,q}
    &=&
    \left(\frac{dg_2}{dy}\right)_q
    =
      \left(D_y g_2\right)_q 
    + \sum_{i=0}^q \left(D_z g_2\right)_{q-i} (Dh_i),
\Label{Dhq-aux3} \\
    \frac{d}{dy}g_{1,q-1-\ell}
    &=&
    \left(\frac{dg_{1}}{dy}\right)_{q-1-\ell}
    =
      \left(D_y g_1\right)_{q-1-\ell} 
    + \sum_{i=0}^{q-1-\ell} \left(D_z g_1\right)_{q-1-\ell-i} (Dh_i).\hskip0.4truein
\Label{Dhq-aux4}
\end{eqnarray}
Substituting
Eqs.~(\ref{Dhq-aux3})
and
(\ref{Dhq-aux4})
into
Eq.~(\ref{Dhq-aux1}),
we obtain
\begin{eqnarray}
      \left(D_y g_2\right)_q 
    + \sum_{i=0}^q \left(D_z g_2\right)_{q-i} Dh_i
    &=&\mbox{}
    \sum_{\ell=0}^{q-1} (D^2 h_\ell) g_{1,q-\ell-1}
    +
    \sum_{\ell=0}^{q-1} 
    (Dh_\ell) (D_y g_1)_{q-1-\ell}
    \nonumber\\
    &+&\mbox{}
    \sum_{\ell=0}^{q-1}\sum_{i=0}^{q-1-\ell} 
    (Dh_\ell) (D_z g_1)_{q-1-\ell-i} (Dh_i).
\Label{Dhq-aux2}
\end{eqnarray}
Separating
the
$i=q$
term in
the sum of
the left member,
changing 
the order of 
summation in
the last sum of
the right member,
and
solving for
$Dh_q$,
we obtain
Eq.~(\ref{Dhq}).
\end{PROOF}

\begin{LEMMA} \Label{DUg-lem}
Let $V$ be a matrix-valued function
of $y$, $z$, and $\eps$
that, together with its first-order derivatives,
is smooth
and
${\cal O}(1)$ as $\eps \downarrow 0$.
If $z=\psif_{(q)} (y, \eps)$ and
\be
    V(\cdot\,,\psif_{(q)}, \eps)
    =
    \sum_{\ell=0}^q
    \eps^\ell V_\ell
    + \mathcal{O}(\eps^{q+1}) ,\quad
    g_1(\cdot\,,\psif_{(q)},\eps)
    =
    \sum_{\ell=0}^q
    \eps^\ell g_{1,_\ell}
    + \mathcal{O}(\eps^{q+1}) ,
\ee
then,
\be
    \frac{dV}{dt} (\cdot\,,\psif_{(q)},\eps)
    =
    \sum_{i=0}^q \eps^{i+1} 
    \sum_{\ell=0}^i\frac{dV_\ell}{dy} g_{1,i-\ell}
    +
    \mathcal{O}(\eps^{q+1}) .
\Label{DVg-exp}
\ee
\end{LEMMA}
\begin{PROOF}
A direct computation gives
\be
  \frac{dV}{dt} = (DV) g
    =
      \eps (D_y V) g_1
    + (D_z V) g_2,
\Label{DUg-aux1}
\ee
where
all the terms
are evaluated at
$(y, \psi_{(q)} (y, \eps),\eps)$.
Since
$\psi_{(q)}$
approximates
the slow manifold
up to and including
$\mathcal{O}(\eps^q)$
terms,
\begin{eqnarray}
    g_1(\cdot\,,\psi_{(q)},\eps)
  = g_1(\cdot\,,h_\eps,\eps) + \mathcal{O}(\eps^{q+1}),
\Label{g1(psi)-vs-g1(h)}\\
    g_2(\cdot\,,\psi_{(q)},\eps)
  = g_2(\cdot\,,h_\eps,\eps) + \mathcal{O}(\eps^{q+1}),
\Label{g2(psi)-vs-g2(h)}
\end{eqnarray}
and also
\be
    D\psi_{(q)} = Dh_\eps + \mathcal{O}(\eps^{q+1}).
\Label{Dpsi-vs-Dh}
\ee
Using
Eqs.~(\ref{inveq}),
(\ref{g1(psi)-vs-g1(h)}),
and
(\ref{Dpsi-vs-Dh}),
we rewrite
Eq.~(\ref{g2(psi)-vs-g2(h)})
as
\be
  g_2(\cdot\,,\psi_{(q)},\eps)
  = \eps (D\psi_{(q)}) g_1(\cdot\,,\psi_{(q)},\eps)
  + \mathcal{O}(\eps^{q+1}).
\Label{inveq-approx(q)-aux}
\ee
Equation~(\ref{inveq-approx(q)-aux})
is an equation for
$\Kf^{(q)}$.
We recast it so
the right member involves
a total derivative
with respect to~$y$,
\be
    (DV) g
    =
      \eps \left(  D_yV + D_zV D\psi_{(q)} \right) 
      g_1
    + \mathcal{O}(\eps^{q+1})
    =
      \eps \frac{dV}{dy} g_1
    + \mathcal{O}(\eps^{q+1})
\Label{DUg-aux3}
\ee
or, expanding in powers of $\eps$,
\begin{eqnarray}
    (DV) g
    &=&
      \sum_{i=0}^q\eps^{i+1} 
      \sum_{\ell=0}^i\left(\frac{dV}{dy}\right)_\ell
      g_{1,i-\ell}
    + 
    \mathcal{O}(\eps^{q+1}).
\Label{DUg-aux4}
\end{eqnarray}
The operations of taking
the total derivative
with respect to~$y$
and expanding
with respect to~$\eps$
commute, so
$(dV/dy)_\ell = dV_\ell/dy$
and Eq.~(\ref{DVg-exp}) follows.
\end{PROOF}
%
%%%%%%%%%%%%%%%%%%%%%%%%%%%%%%%%%%%%%%%%%%%%
%

\newpage
\noindent
Corresponding author:

\noindent
Hans G.\ Kaper \\
Division of Mathematical Sciences \\
National Science Foundation \\
4201 Wilson Boulevard, Suite 1025 \\
Arlington, VA 22230

\noindent
Authors' e-mail addresses:

\noindent
\texttt{azagaris@math.bu.edu} \\
\texttt{kaper@mcs.anl.gov; hkaper@nsf.gov} \\
\texttt{tasso@math.bu.edu} \\

\vfill
\begin{flushright}
\scriptsize
\framebox{\parbox{2.4in}{
The submitted manuscript has been created
by the University of Chicago as Operator of
Argonne National Laboratory ("Argonne")
under Contract No.\ W-31-109-ENG-38
with the U.S.\ Department of Energy.
The U.S.\ Government retains for itself,
and others acting on its behalf,
a paid-up, nonexclusive, irrevocable
worldwide license in said article
to reproduce, prepare derivative works,
distribute copies to the public,
and perform publicly and display publicly,
by or on behalf of the Government.}}
\normalsize
\end{flushright}

\end{document}